\documentclass[format=acmsmall,screen]{acmart}
\usepackage{acm-ec-23}
\usepackage{booktabs} 
\usepackage[ruled]{algorithm2e} 

\SetAlFnt{\small}
\SetAlCapFnt{\small}
\SetAlCapNameFnt{\small}
\SetAlCapHSkip{0pt}
\IncMargin{-\parindent}

\setcitestyle{authoryear}

\newcommand{\cA}{{\mathcal{A}}}
\newcommand{\cU}{{\mathcal{U}}}
\newcommand{\cG}{{\mathcal{G}}}
\newcommand{\cV}{{\mathcal{V}}}

\newcommand{\cS}{{\mathcal{S}}}

\newcommand{\cL}{{\mathcal{L}}}
\newcommand{\cX}{{\mathcal{X}}}

\newcommand{\cF}{{\mathcal{F}}}

\newcommand{\cT}{{\mathcal{T}}}

\newcommand{\cC}{{\mathcal{C}}}
\newcommand{\cH}{{\mathcal{H}}}

\newcommand{\cO}{{\mathcal{O}}}
\newcommand{\cE}{{\mathcal{E}}}

\newcommand{\bbe}{{\textbf{e}}}
\newcommand{\bI}{{\textbf{I}}}
\newcommand{\bJ}{{\textbf{J}}}
\newcommand{\bM}{{\textbf{M}}}

\newcommand{\bB}{\textbf{B}}

\newcommand{\bR}{\textbf{R}}

\newcommand{\bY}{\textbf{Y}}
\newcommand{\bx}{\textbf{x}}
\newcommand{\by}{\textbf{y}}

\newcommand{\bq}{\textbf{q}}
\newcommand{\bu}{\textbf{u}}

\newcommand{\bH}{\textbf{H}}
\newcommand{\bZ}{\textbf{Z}}

\newcommand{\bp}{\textbf{p}}
\newcommand{\ba}{\textbf{a}}
\newcommand{\bv}{\textbf{v}}

\newcommand{\bb}{\textbf{b}}

\newcommand{\bbeta}{{\pmb{\eta}}}
\newcommand{\bgamma}{{\pmb{\gamma}}}

\newcommand{\bbR}{\mathbb{R}}

\newcommand{\bbI}{\mathbb{I}}
\newcommand{\bbN}{\mathbb{N}}

\usepackage{mathtools}

\newcommand\myeq{\stackrel{\mathclap{\tiny\mbox{def}}}{=}}

\newtheorem{theorem}{Theorem}
\newtheorem{lemma}{Lemma}
\newtheorem{remark}{Remark}

\usepackage{booktabs}

\usepackage[normalem]{ulem}

\newif\ifnotes\notestrue
\def\remove#1{\sout{#1}}
\def\remove#1{}

\def\mtien#1{{\color{black}{#1}}}

\def\htien#1{}

\usepackage[ruled]{algorithm2e} 
\usepackage{color}

\usepackage{tikz}

\title[Tackling Stackelberg Network Interdiction under Boundedly Rational Adversary]{Tackling Stackelberg Network Interdiction against a Boundedly Rational Adversary}

\author{Tien Mai}
\affiliation{%
	\institution{Singapore Management University}
	\city{Singapore}\country{Singapore}
}
\author{Avinandan Bose}
\affiliation{%
	\institution{University of Washington}
	\city{Washington}
	\country{United States}
}
\author{Arunesh Sinha}
\affiliation{%
	\institution{Rutgers University}
	\city{}
	\country{United States}
}
\author{Thanh H. Nguyen}
\affiliation{%
	\institution{University of Oregon}
	\city{}
	\country{United States}
}
\author{}

\begin{abstract}
 This work studies Stackelberg network interdiction games --- an important class of games in which a defender first allocates (randomized) defense resources to a set of critical nodes on a graph while an adversary chooses its path to attack these nodes accordingly. 
    We consider a boundedly rational adversary in which the adversary's response model is based on a dynamic form of classic logit-based discrete choice models. We show that the problem of finding an optimal interdiction strategy for the defender in the  rational setting is NP-hard. The resulting optimization is in fact non-convex and additionally, involves complex terms that sum over exponentially many paths. We tackle these computational challenges by presenting new efficient approximation algorithms with bounded solution guarantees. First, we address the exponentially-many-path challenge by proposing a polynomial-time dynamic programming-based formulation. We then show that the gradient of the non-convex objective can also be computed in polynomial time, which allows us to use a gradient-based method to solve the problem efficiently. 
    Second, we identify a restricted problem that is convex and hence gradient-based methods find the global optimal solution for this restricted problem. We further identify mild conditions under which this restricted problem provides a bounded approximation for the original problem.
\end{abstract}

\begin{document}

\begin{titlepage}

\maketitle

\end{titlepage}

\section{Introduction}

Network interdiction is a well-studied topic in Artificial Intelligence. 
There are many practical problems~\cite{smith2020survey}, such as in cyber systems and illicit supply networks, that can be modelled as a network interdiction problem. In literature, many variations in models of network interdiction exist, and consequentially, a variety of techniques have been used for solving different types of these problems. Our work focuses on a particular type in which there is a set of critical nodes $\cL$ to protect within a larger network of $\cS$ nodes.
We employ a popular network interdiction  model~\cite{fulkerson1977maximizing,israeli2002shortest},  
where the interdictor (defender) uses a randomized allocation of limited defense resources for the critical nodes in $\cL$. The adversary traverses the graphs starting from an origin $s_o$ and reaching a destination $s_d$. There is an interaction with the defender only if the adversary crosses any node in $\cL$. The interaction is modelled using a leader-follower (Stackelberg) game setting where the defender first allocates resources in a randomized fashion and then the adversary chooses its path accordingly. 

We model the adversary behavior using a dynamic Quantal Response model (an instance of well-known dynamic discrete choice models \citep{RUST1987,AGUIRREGABIRIA2010}). 
In this model, the adversary, currently at node $s$, chooses to visit a next node $s'$ with a probability proportional to the quantity $\exp\left(\frac{v}{\mu}\right)$ where $v$ is the adversary's immediate utility of choosing $s'$ among neighbors of $s$ and $\mu$ is a constant. In the game context, the utility $v$ depends on the defender resource allocation denoted by $\bx$ as $v(\cdot|\bx)$. In a seminal result \citep{FosgFrejKarl13,MaiFosFre15}, it was shown that, under some specific settings, such local Quantal Response choices directly correspond to a Quantal Response choice over paths in the network graph. That is, the probability of choosing one path $\tau$ from the origin to the destination is proportional to $\exp\left(\frac{U(\tau)}{\mu}\right)$ where $U(\tau)$ is the sum of utilities along the path. We use this adversary response model to form an optimization for computing an optimal interdiction strategy of the defender. To the best of our knowledge, existing  network interdiction models assume perfectly rational adversaries and make use of some linear programming techniques to handle
\citep{smith2009new,smith2020survey}. We are the first to explore the DDC framework to model bounded rational adversaries and formulate the defender's problem as nonlinear optimization ones, opening the possibility of solving the network interdiction problems via nonlinear optimization techniques.

While the closed form result of the adversary response is mathematically interesting, it presents computational challenges as the computation of any such probability involves reasoning about exponentially many paths from origin to destination. 
In this paper, we show that it is NP-hard to find an optimal interdiction strategy for the defender in this game setting of a boundedly rational adversary. Therefore, We address the challenge of solving such complex non-convex optimization problem for the defender with the following new efficient approximation algorithms.

First, we propose an efficient dynamic programming method to compute the objective (defender expected utility) as well as gradient of objective of the optimization even though these terms involve summing over exponentially many paths. This is accomplished by exploiting recursive relationships among adversary utility-related terms across different paths that involves in the defender's optimal interdiction strategy computation. By employing dynamic programming, we can follow a gradient descent approach that is computationally efficient at each step to optimize the defender strategy. 

Second, while the above proposed method is computationally efficient, it does not guarantee global optimality due to the non-convexity nature of the defender's optimization problem. Therefore, we identify a restricted problem in which the adversary can visit only one (any one) node in $\cL$ and show that the optimization is convex with such a restriction. Therefore, this restricted problem can be solved optimally in a tractable manner using the efficient gradient descent from the first contribution. 
We further identify conditions on two problem specific terms $\beta_1, \beta_2$ (Equation~\ref{eq:beta}) such that if $\beta_1, \beta_2$ are small, the solution to the above restricted problem provides close approximation guarantees for the original unrestricted problem. 

\noindent
\textbf{Notation:}
Boldface characters represent matrices or vectors or sets, and $a_i$ denotes the $i$-th element of $\ba$ if $\ba$ is indexable. We use $[m]$, for any $m\in \mathbb{N}$, to denote the set $\{1,\ldots,m\}$. 
\section{Related Work}
\noindent \textbf{Dynamic discrete choice models.} We employ the dynamic discrete choice (DDC) framework to model the adversary bounded rational behavior. From the seminal work of \citep{RUST1987}, DDC models have been widely studied and used to analyze  sequential looking-forward choice behaviors and have various applications, e.g., on fertility and child mortality~\citep{Wolpin1984}, on job matching and occupational choice~\cite{Miller1984}, on bus engine replacement~\cite{RUST1987}, and on route choice analysis~\cite{FosgFrejKarl13,MaiFosFre15}. Among existing DDC models, the logit-based DDC has been popular due to its closed-form formulation~\citep{RUST1987}. This model can be viewed as a dynamic version of the well-known multinomial logit (or Quantal Response) model \citep{McFa81,Trai03}. In transportation modeling, or specially route choice analysis, the logit-based DDC model was utilized under an undiscounted infinite horizon Markov  Decision Process to develop models to predict people's bounded rational path-choice behavior~\cite{FosgFrejKarl13,MaiFosFre15}. As highlighted in \cite{zimmermann2020tutorial}, such a route choice  model presents synergies with the
\textit{stochastic shortest path}  problem \cite{bertsekas1991analysis}. 

\noindent \textbf{Network interdiction.} Our work is closely related to the well-studied shortest path interdiction problem~\citep{fulkerson1977maximizing,israeli2002shortest} and can be viewed as its bounded rational version. The shortest path interdiction and other network interdiction problems with perfectly rational adversaries are generally NP-hard and have strong connections with the areas of bi-level optimization \citep{dempe2015bilevel}  and robust optimization \citep{ben2002robust}. We refer the reader to  \cite{smith2020survey} for a comprehensive review. As we mentioned previously, existing network interdiction models consider perfectly rational adversaries
\citep{smith2009new,smith2020survey}. We, on the other hand, explore the DDC framework to model bounded rational adversaries, resulting in a significantly more challenging defender problem as it involves complex nonlinear optimization.  
Besides, there are other variant models where the problem data is not perfectly known to players \citep{cormican1998stochastic}, or where the players repeatedly make their actions alternatively \citep{sefair2016dynamic}, or where online learning is involved \cite{borrero2016sequential}. 
These models provide promising next steps for future work.

\noindent \textbf{Network security games and others.} Our work also relates to static Stackelberg security game models with Quantal Response adversaries \citep{yang2011improving,yang2012computing,1501625,mai2022choices,vcerny2021dinkelbach,milec2020complexity}. In dynamic models named as network security games~\cite{jain2011double}, the set-up is different from our work as in this work the \emph{rational} adversary aims to reach a target and stop, whereas in our work the \emph{boundedly rational} adversary can attack multiple targets. Other related works along this line only consider \emph{zero-sum} network security game setting; they attempt to develop scalable game-theoretic algorithms using techniques in either game theory or machine learning~\cite{xue2021solving,xue2022nsgzero}. A Quantal Response type relaxation for network security game was also studied, where the focus in on smart predict and optimize~\cite{wang2020scalable}, however, the optimize part is done using standard non-linear solver such sequential quadratic program with no guarantees. 

In addition, there are other related game models to our work in the sense that players act in a graph-based environment, including pursuit-evasion games and security patrol games~\cite{zhang2019optimal,basilico2009leader,basilico2017adversarial}. However, these existing works do not consider the attacker's bounded rationality. Strategy spaces are also characterized differently in these works where real-time information or alarm signals., etc are incorporated into defining strategies of players.  

\section{Problem Formulation}
Our network interdiction problem is a leader-follower game with a single adversary. The game is played on a network (graph). We formulate the problem as a two-agent \remove{Markov Decision Process} \mtien{network interdiction game} \remove{ where the nodes are states, hence} where the set of nodes \remove{(or state space)} is given as $\cS$. The follower (adversary) takes a path through this network, which is sampled from a distribution as described below.
Let $\cL$ be the set of critical \mtien{nodes} (i.e., subset of nodes in the network) that the defender can interfere or alter.
From the leader's (defender's) viewpoint, the aim is to 
assign $M$ resources to nodes $s \in \cL$; each such assignment is a defender pure strategy. Further, nodes and resources are of certain kinds such that nodes of a given kind can only be protected by resources of that same kind. Let there be $K$ kinds of nodes and corresponding resources. Let the number of resources of each kind $k$ be $M_k$, hence $\sum_{k \in [K]} M_k = M$. Also, let $\{\cL_k\}_{k \in [K]}$ be a partition of the set of nodes $\cL$ by the kind of the nodes.

A mixed strategy is a randomized allocation resulting in marginal probability $x_s$ of covering node $s$ (with $\bx=\{x_s,\; s\in \cL\}$), with the restriction that \mtien{$\sum_{s \in \cL_k} x_s \leq M_k$} for all $k \in [K]$,
which then impacts the adversary's path choice probabilities.
  Given a \remove{state} node $s \in \cS$, if the adversary crosses this \remove{state} node, then the defender gets a node-specific reward $r^l(s, x_s)$. The defender's expected reward for following the interdiction strategy $\bx$ can be computed as follows:
\[
\cF^l(\bx) =  \sum_{s\in \cL}r^l(s, x_s) P^f_s(\bx),
\]
where $P^f_s(\bx)$ is the probability that the follower  (adversary) crosses a \remove{state} \mtien{node} $s$, computed as:
\[
P^f_s(\bx)  = \sum_{\tau \in  \Omega, \tau \ni s} \prod_{(s_t,s_{t+1}) \in \tau}\pi^f(s_{t+1}|s_{t},\bx),
\]
where $\Omega$ is the set of all possible \remove{trajectories} {paths} and $\pi^f(s_{t+1}|s_t,\bx)$ is the randomized path selection policy of the adversary. Note that the defender's interdiction action $\bx$ is fixed and known to the adversary when the adversary is playing.
The equation for $P^f_s(\bx)$ is essentially the sum of probabilities that the adversary will actually move along each trajectory $\tau$ that crosses \remove{the state} {node} $s$  $(\tau \ni s)$. Therefore, the adversary's policy has a Markovian nature.
The optimization problem to find an optimal defender interdiction strategy can be formulated as follows:
\begin{align}
     \max_{\bx } &\quad\quad \cF^l(\bx) =  \sum_{s\in \cL}r^l(s, x_s) P^f_s(\bx) \label{prob:main-homo} \tag{\sf OPT}\\
    \mbox{subject to} & \quad\quad \sum_{s \in \cL_k} x_s \leq M_k, \; \forall k \in [K]  \text{ and }x_s \in[L^x,U^x],\;\forall s\in \cL, \nonumber
\end{align}
which is generally non-convex. Here, $[L^x, U^x]$ represent the required lower bound and upper bound on the coverage probability for each \remove{state} {node} in the critical set $\mathcal{L}$.

\paragraph{\textbf{Boundedly rational adversary behavior in network interdiction games.}} We model the adversary's bounded rational behavior using the dynamic discrete choice framework {(and specifically the logit-based recursive path choice model \citep{FosgFrejKarl13})}. To describe formally, let us consider a network $(\cS,\cA)$ where $\cS$ is a set of nodes $\cS = \{1,2,\ldots,|\cS|\}$, and $\cA$ is a set of arcs. 
\remove{deterministic discrete Markov Decision Process (MDP) defined by a tuple $(\cS,\cA,\bq, \bu, s_o)$, where $\cS$ is a set of states (nodes) $\cS = \{1,2,\ldots,|\cS|\}$, $\cA$ is a finite set of actions.} \remove{In addition,
 $\bq$ is a matrix of transition probabilities given by individual elements
$q:\cS\times \cA\times\cS \rightarrow [0,1]$, i.e., $q(s_{t+1}|a_t,s_t)$ is the probability of moving to state $s_{t+1}\in\cS$ from $s_t\in \cS$ by performing action  $a_t\in \cA$ at time step $t$.} 
Moreover, let $\bv $ be a matrix of utilities of the adversary, i.e., each $v(s; \bx)$ is the immediate utility of visiting node  $s\in \cS$,  given leader strategy $\bx$. The origin $s_o \in \cS$ is a given starting node. In our problem, we also assume the existence of a sink (or destination) node $s_d \in \cS$ that the adversary ultimately reaches. {This can be a physical node of the network, or a dummy one representing the final state of the adversary.}

Under the logit-based dynamic discrete choice framework (or dynamic Quantal Response), we assume that the adversary responds in a bounded rational manner.
Let $U(\tau | s_0, \bx) = \sum_{s\in\tau }v(s;\bx)$ be the deterministic long-term utility when starting in $s_0$; if $s_0 = s_o$, then we simply write $U(\tau | \bx)$. 
A known property in this setting is that the bounded rational adversary chooses a policy that is equivalent to a static multinomial logit (MNL) discrete choice model (or logit Quantal Response model) over all possible paths~\cite{FosgFrejKarl13}.
More precisely, given $\bx$,  
the probability that the adversary follows a \remove{trajectory} {path} $\tau$  can be computed as follows~\citep{FosgFrejKarl13}:
\begin{equation}
\label{eq:tau-prob}
P(\tau) =  \frac{\exp \left(\frac{U(\tau;\bx)}{\mu}\right)}{Z}, \text{ where }Z = \sum_{\tau\in \Omega} \exp \left(\frac{U(\tau;\bx)}{\mu}\right), 
\end{equation}
given $\Omega$ is the set of all possible {paths} and $\mu$ is the parameter which governs the follower's rationality. \remove{Note that, given deterministic transitions, $\tau$ can be thought of a path $\{s_0, s_1, s_2 ,\ldots\}$ in the state space graph with consecutive states forming the edge $(s,s')$ with weight $v(s;\bx)$. }
Thus, we can view the logit-based dynamic discrete choice formulation as a soft version of the shortest weighted path problem from the  source $s_o$ to { destination $s_d$.} \remove{on the state space graph.} Given the adversary behavior model, the adversary's expected utility can be computed as follows:
\begin{align*}
    \cE^f(\bx) = \sum_{\tau\in \Omega}P(\tau)U(\tau; \bx)
\end{align*}
which is the expectation over all paths. We present the following proposition which shows that the adversary's expected utility approaches the best  accumulated utility (smallest path weight) as $\mu$ tends to zero (we drop the fixed $\bx$ for simplicity).

\begin{proposition}\label{prop.1}
Let $\tau^* = \text{argmax}_{\tau \in \Omega}U(\tau)$ (i.e., the best {path} which gives the highest adversary utility) and $L^* = \max_{\tau}|U(\tau)|$. In addition, let $\Omega^* = \{\tau;\; U(\tau) = U(\tau^*)\}$ (i.e., the set of all \remove{trajectories} {paths} with the same highest utility) and $\alpha = U(\tau^*) - \max_{\tau \in \Omega\backslash\Omega^*} U(\tau)$ (i.e., the utility loss of the adversary if he selects the second best {path} instead of the best one). Then we obtain:
\[
 {|\cE^f -  U(\tau^*)| \leq \frac{L^*+1}{1+\frac{|\Omega^*|}{|\Omega\backslash\Omega^*|}\exp\left(\frac{\alpha}{\mu}\right)}}.
\]
As a result, $\lim_{\mu\rightarrow 0 }\cE^f =  U(\tau^*)$.\footnote{All proofs, if not presented, are included in the appendix.}
\end{proposition}

Based on the above proposition $\mu=0$ results in an utility maximizing (rational) adversary. Our Theorem~\ref{thm.NPhard} shows that \eqref{prob:main-homo} is NP-Hard for a rational adversary and hence we focus on approximations (with bounded solution guarantees) in the rest of the paper.
\begin{theorem}\label{thm.NPhard}
    The defender interdiction problem \eqref{prob:main-homo} is NP-Hard for $\mu = 0$.
\end{theorem}
\begin{proof}
    We do a reduction from exact 3-Cover problem,  where given $m$ items $\{1, \ldots, m\}$ and a collection of $n$ subsets $\{S_1, \ldots, S_n\}$ with $S_v \subset \{1, \ldots, m\}$ each of size 3, i.e., $|S_v| = 3$ for $v \in \{1, \ldots, n\}$, the decision problem is whether there is a cover that contains each item exactly once. This is a known NP-Hard problem. Also, clearly any valid cover must be of $m/3$ number of subsets.
\paragraph{Game instance construction given the 3-Cover problem}
Given an exact 3-Cover problem, we form an instance of our network security game as follows: 
The critical nodes can be one of $m+1$ types, among which the first $m$ types are labeled $1, \ldots, m$ and the last type is labeled \emph{red}. In addition,
there are a total of $n - m/3 + m$ defender resources. Among these resources, $n-m/3$ resources, denoted by $R_1, \ldots, R_{n-m/3}$, can defend nodes of type red (we call these the red resources). The remaining $m$ resources are denoted by $r_1, \ldots, r_m$, where resource $r_j$ can defend a node of type $j$. 

We form $n$ sub-graphs --- each sub-graph corresponds to a subset $S_v = \{i,j,k\}$ shown below in the picture. The sub-graph for $S_v$ has one critical node of type red and three critical nodes labeled $(v,i), (v,j), (v,k)$ of types $i,j,k$ respectively. There is an initial non-critical node $s_0$ and an end non-critical node $s_e$ in the sub-graph. A direct edge also connects $s_0$ to $s_e$, called a dummy edge. We can join all the $n$ sub-graphs by making a source node and connect the source node to all $s_0$ for every sub-graph and a sink node and connect all $s_e$ of each sub-graph to the sink node (see below).

\begin{center}

\tikzset{every picture/.style={line width=0.75pt}} 

\begin{tikzpicture}[x=0.75pt,y=0.75pt,yscale=-0.8,xscale=0.80]

\draw   (122,81) .. controls (122,73.27) and (128.27,67) .. (136,67) .. controls (143.73,67) and (150,73.27) .. (150,81) .. controls (150,88.73) and (143.73,95) .. (136,95) .. controls (128.27,95) and (122,88.73) .. (122,81) -- cycle ;
\draw  [fill={rgb, 255:red, 237; green, 72; blue, 72 }  ,fill opacity=1 ] (213,78) .. controls (213,70.27) and (219.27,64) .. (227,64) .. controls (234.73,64) and (241,70.27) .. (241,78) .. controls (241,85.73) and (234.73,92) .. (227,92) .. controls (219.27,92) and (213,85.73) .. (213,78) -- cycle ;
\draw   (392,80) .. controls (392,72.27) and (398.27,66) .. (406,66) .. controls (413.73,66) and (420,72.27) .. (420,80) .. controls (420,87.73) and (413.73,94) .. (406,94) .. controls (398.27,94) and (392,87.73) .. (392,80) -- cycle ;
\draw   (294,130.5) .. controls (294,121.94) and (300.94,115) .. (309.5,115) .. controls (318.06,115) and (325,121.94) .. (325,130.5) .. controls (325,139.06) and (318.06,146) .. (309.5,146) .. controls (300.94,146) and (294,139.06) .. (294,130.5) -- cycle ;
\draw   (291,42.5) .. controls (291,33.39) and (298.39,26) .. (307.5,26) .. controls (316.61,26) and (324,33.39) .. (324,42.5) .. controls (324,51.61) and (316.61,59) .. (307.5,59) .. controls (298.39,59) and (291,51.61) .. (291,42.5) -- cycle ;
\draw   (292,86) .. controls (292,77.16) and (299.16,70) .. (308,70) .. controls (316.84,70) and (324,77.16) .. (324,86) .. controls (324,94.84) and (316.84,102) .. (308,102) .. controls (299.16,102) and (292,94.84) .. (292,86) -- cycle ;
\draw    (233.5,69) -- (289.18,43.34) ;
\draw [shift={(291,42.5)}, rotate = 155.26] [color={rgb, 255:red, 0; green, 0; blue, 0 }  ][line width=0.75]    (10.93,-3.29) .. controls (6.95,-1.4) and (3.31,-0.3) .. (0,0) .. controls (3.31,0.3) and (6.95,1.4) .. (10.93,3.29)   ;
\draw    (150,81) -- (211,78.1) ;
\draw [shift={(213,78)}, rotate = 177.27] [color={rgb, 255:red, 0; green, 0; blue, 0 }  ][line width=0.75]    (10.93,-3.29) .. controls (6.95,-1.4) and (3.31,-0.3) .. (0,0) .. controls (3.31,0.3) and (6.95,1.4) .. (10.93,3.29)   ;
\draw    (238,78) -- (291.01,83.78) ;
\draw [shift={(293,84)}, rotate = 186.23] [color={rgb, 255:red, 0; green, 0; blue, 0 }  ][line width=0.75]    (10.93,-3.29) .. controls (6.95,-1.4) and (3.31,-0.3) .. (0,0) .. controls (3.31,0.3) and (6.95,1.4) .. (10.93,3.29)   ;
\draw    (235.5,86) -- (292.41,129.29) ;
\draw [shift={(294,130.5)}, rotate = 217.26] [color={rgb, 255:red, 0; green, 0; blue, 0 }  ][line width=0.75]    (10.93,-3.29) .. controls (6.95,-1.4) and (3.31,-0.3) .. (0,0) .. controls (3.31,0.3) and (6.95,1.4) .. (10.93,3.29)   ;
\draw    (324,42.5) -- (391.68,73.17) ;
\draw [shift={(393.5,74)}, rotate = 204.38] [color={rgb, 255:red, 0; green, 0; blue, 0 }  ][line width=0.75]    (10.93,-3.29) .. controls (6.95,-1.4) and (3.31,-0.3) .. (0,0) .. controls (3.31,0.3) and (6.95,1.4) .. (10.93,3.29)   ;
\draw    (324,84) -- (390,80.12) ;
\draw [shift={(392,80)}, rotate = 176.63] [color={rgb, 255:red, 0; green, 0; blue, 0 }  ][line width=0.75]    (10.93,-3.29) .. controls (6.95,-1.4) and (3.31,-0.3) .. (0,0) .. controls (3.31,0.3) and (6.95,1.4) .. (10.93,3.29)   ;
\draw    (325,130.5) -- (394.76,90.99) ;
\draw [shift={(396.5,90)}, rotate = 150.47] [color={rgb, 255:red, 0; green, 0; blue, 0 }  ][line width=0.75]    (10.93,-3.29) .. controls (6.95,-1.4) and (3.31,-0.3) .. (0,0) .. controls (3.31,0.3) and (6.95,1.4) .. (10.93,3.29)   ;
\draw    (143.5,93) .. controls (312.79,226.65) and (364.96,126.06) .. (404.8,94.92) ;
\draw [shift={(406,94)}, rotate = 143.13] [color={rgb, 255:red, 0; green, 0; blue, 0 }  ][line width=0.75]    (10.93,-3.29) .. controls (6.95,-1.4) and (3.31,-0.3) .. (0,0) .. controls (3.31,0.3) and (6.95,1.4) .. (10.93,3.29)   ;
\draw   (38,196) .. controls (38,182.19) and (49.19,171) .. (63,171) .. controls (76.81,171) and (88,182.19) .. (88,196) .. controls (88,209.81) and (76.81,221) .. (63,221) .. controls (49.19,221) and (38,209.81) .. (38,196) -- cycle ;
\draw   (493,183) .. controls (493,169.19) and (504.19,158) .. (518,158) .. controls (531.81,158) and (543,169.19) .. (543,183) .. controls (543,196.81) and (531.81,208) .. (518,208) .. controls (504.19,208) and (493,196.81) .. (493,183) -- cycle ;
\draw  [dash pattern={on 4.5pt off 4.5pt}] (116,20) -- (454.5,20) -- (454.5,170) -- (116,170) -- cycle ;
\draw    (71.5,172) -- (125.4,90.67) ;
\draw [shift={(126.5,89)}, rotate = 123.53] [color={rgb, 255:red, 0; green, 0; blue, 0 }  ][line width=0.75]    (10.93,-3.29) .. controls (6.95,-1.4) and (3.31,-0.3) .. (0,0) .. controls (3.31,0.3) and (6.95,1.4) .. (10.93,3.29)   ;
\draw    (417.5,85) -- (502.02,161.66) ;
\draw [shift={(503.5,163)}, rotate = 222.21] [color={rgb, 255:red, 0; green, 0; blue, 0 }  ][line width=0.75]    (10.93,-3.29) .. controls (6.95,-1.4) and (3.31,-0.3) .. (0,0) .. controls (3.31,0.3) and (6.95,1.4) .. (10.93,3.29)   ;
\draw  [dash pattern={on 4.5pt off 4.5pt}] (117,183) -- (455.5,183) -- (455.5,258) -- (117,258) -- cycle ;
\draw    (451.5,216) -- (495.64,198.73) ;
\draw [shift={(497.5,198)}, rotate = 158.63] [color={rgb, 255:red, 0; green, 0; blue, 0 }  ][line width=0.75]    (10.93,-3.29) .. controls (6.95,-1.4) and (3.31,-0.3) .. (0,0) .. controls (3.31,0.3) and (6.95,1.4) .. (10.93,3.29)   ;
\draw  [dash pattern={on 4.5pt off 4.5pt}]  (289.5,265) -- (289.92,295.98) -- (290.5,339) ;
\draw    (80.5,214) -- (122.52,219.73) ;
\draw [shift={(124.5,220)}, rotate = 187.77] [color={rgb, 255:red, 0; green, 0; blue, 0 }  ][line width=0.75]    (10.93,-3.29) .. controls (6.95,-1.4) and (3.31,-0.3) .. (0,0) .. controls (3.31,0.3) and (6.95,1.4) .. (10.93,3.29)   ;
\draw    (68,219) -- (127.55,329.24) ;
\draw [shift={(128.5,331)}, rotate = 241.62] [color={rgb, 255:red, 0; green, 0; blue, 0 }  ][line width=0.75]    (10.93,-3.29) .. controls (6.95,-1.4) and (3.31,-0.3) .. (0,0) .. controls (3.31,0.3) and (6.95,1.4) .. (10.93,3.29)   ;
\draw    (446.5,322) -- (507.55,208.76) ;
\draw [shift={(508.5,207)}, rotate = 118.33] [color={rgb, 255:red, 0; green, 0; blue, 0 }  ][line width=0.75]    (10.93,-3.29) .. controls (6.95,-1.4) and (3.31,-0.3) .. (0,0) .. controls (3.31,0.3) and (6.95,1.4) .. (10.93,3.29)   ;
\draw  [dash pattern={on 4.5pt off 4.5pt}]  (106.5,220) .. controls (113.5,231) and (100.5,251) .. (83.5,248) ;
\draw  [dash pattern={on 4.5pt off 4.5pt}]  (492.5,233) .. controls (478.5,239) and (473.82,238.07) .. (468.63,230.24) .. controls (463.45,222.4) and (465.65,214.26) .. (466.5,210) ;

\draw (130,70.4) node [anchor=north west][inner sep=0.75pt]    {$s_{0}$};
\draw (295,34.4) node [anchor=north west][inner sep=0.75pt]    {$v,i$};
\draw (296,77.4) node [anchor=north west][inner sep=0.75pt]    {$v,j$};
\draw (297,120.4) node [anchor=north west][inner sep=0.75pt]    {$v,k$};
\draw (399,70.4) node [anchor=north west][inner sep=0.75pt]    {$s_{e}$};
\draw (39,186) node [anchor=north west][inner sep=0.75pt]   [align=left] {Source};
\draw (502,174) node [anchor=north west][inner sep=0.75pt]   [align=left] {Sink};

\end{tikzpicture}
\end{center}

For any node, $p$ denotes the probability the defender protects that node.
We set the payoff of the adversary for the critical nodes as $u_a(p) = -50 p + 50 (1-p) - K p \log p = 50 - 100p - K p \log p$ where the constant $K =\frac{100}{\log\left(\frac{n-m/3+1}{n-m/3 + 0.5}\right)} > 0$. The attacker payoff of skipping all critical nodes via the bottom edge in each sub-graph is 0. 
On the other hand, the defender's payoff, when an adversary visits a critical node is set to $u_d(p) = -100 (1-p) - \epsilon$ for some small $\epsilon >0$. Thus, the defender's payoff is always strictly negative if the attacker crosses any critical node. When the adversary does not visit any critical node, the defender payoff is $0$. This means that any defender optimal strategy gives the defender a maximum payoff of at most $0$.
\paragraph{{Problem reduction.}}

We claim that there exists an exact 3-Cover if and only if the defender's optimal expected payoff for the corresponding network interdiction game is $0$. 

First, assume there is an exact 3-Cover, then we show that the following strategy provides a payoff of $0$ to the defender: (i) the defender allocates the $n-m/3$ red resources $R_1, \ldots, R_{n-m/3}$ to the red nodes in the $n-m/3$ sub-graphs corresponding to non-cover subsets; and (ii) for the $m/3$ sub-graphs corresponding to subsets in cover, the defender allocates the $m$ resources $r_1, \ldots, r_m$ to the three nodes in each sub-graph. This ensures that either the red node or the three nodes $(v,i), (v,j),(v,k)$ are completely protected in every sub-graph $S_v$. Given this strategy of the defender, we note that in any non-dummy path, there are exactly two critical nodes --- one node will be protected by the defender with a probability of $1$ and the other critical node is protected with a probability of zero. As a result, if the adversary chooses this path, the adversary will obtain a total expected payoff over these two nodes as equal to $-50 + 50 = 0$ (given $p \log p = 0$ when $p$ is either $0$ or $1$). Breaking ties in favor of defender~\cite{leitmann1978generalized}, the adversary will choose one of the bottom dummy edges, providing an expected payoff of zero for the defender. As we discussed previously, the maximum payoff the defender can achieve is at most $0$. Therefore, the above strategy is an optimal strategy for the defender that leads to the maximum defender payoff of $0$.

Next, assume that the defender gets an optimal expected payoff of $0$ and the equilibrium strategy is a vector of probability values for each critical node $\bp^*$. As the expected payoff is $0$, the adversary must have chosen one of the bottom edges. Let us analyze one path through the critical nodes that has $p^*_{v,r}$ on red node and $p_{v,i}^*$ on the other node of type $i$. The adversary payoff for choosing this path is $100 - 100(p^*_{v,r} + p_{v,i}^*) - K p^*_{v,r} \log p^*_{v,r} - K p_{v,i}^* \log p_{v,i}^*$. This adversary payoff for this path must be $\leq 0$ (since adversary chooses the bottom edge). Or by rearranging and dividing by 100,
\begin{align} \label{eq:hardness1}
p^*_{v,r} (1 + (K/100) \log p^*_{v,r}) + p_{v,i}^* (1 + (K/100) \log p_{v,i}^*) \geq 1
\end{align}
Based on Eq.~\ref{eq:hardness1}, we are going to show that the red nodes are uncovered (i.e., $p^*_{v,r}=0$) for exactly $m/3$ sub-graphs. First, since we have $n-m/3$ defender resources that can cover red nodes and we have $n$ red nodes in total (i.e., one red node for each of $n$ sub-graphs), the number of red nodes are uncovered (i.e., $p^*_{v,r}=0$) is $\leq m/3$. Furthermore, we will show by contradiction that the other situation of $p^*_{v,r} = 0$ for $< m/3$ sub-graphs cannot occur. To prove by contradiction, assume that $p^*_{v,r} = 0$ for $< m/3$ sub-graphs or in other words $p^*_{v,r} > 0$ for $\geq n- m/3 + 1$ graphs. In the following, we prove that this assumption leads to a contradiction and hence this assumption cannot hold.

We observe that $p_{v,i}^* (1 + (K/100) \log p_{v,i}^*) \leq 1$. Thus, from Eq.~\ref{eq:hardness1}, we obtain: 
\begin{align} \label{eq:hardness2}
p^*_{v,r} (1 + (K/100) \log p^*_r)  \geq 0
\end{align}

Let $P \subset \{1,\ldots, n\}$ be the subset for which $p^*_{v,r} > 0$ for $v \in P$ (and by our assumption $|P| \geq n- m/3 + 1$). 
For any such $p^*_{v,r} > 0$ with $v \in P$, from Eq.~\ref{eq:hardness2}, we must have $1 + (K/100) \log p^*_{v,r} \geq 0$, which when rearranged gives $p^*_{v,r} \geq \exp(-100/K)$. Summing over all $v \in P$, we get $$\sum_{v \in P} p^*_{v,r} \geq |P| \exp(-100/K)$$ 
The LHS above is $\sum_{v \in P} p^*_{v,r} \leq n - m/3$ (as $p^*_{v,r} = 0$ for $v \notin P$ and there are only $n-m/3$ red resources). With $K =\frac{100}{\log\left(\frac{n-m/3+1}{n-m/3 + 0.5}\right)} $ and the fact that $|P| \geq n- m/3 + 1$, we get the RHS is $\geq n-m/3 + 0.5$. This is a contradiction. Thus, the assumption that $p^*_{v,r} = 0$ for $< m/3$ sub-graphs does not hold.  

Therefore, the red nodes are uncovered (i.e., $p^*_{v,r}=0$) for exactly $m/3$ sub-graphs. We note that if $p^*_{v,r} = 0$, then $p_{v,i}^* = 1$ in order to satisfy Eq.~\ref{eq:hardness1}. And with the same reasoning for other non-dummy paths in the same sub-graph $S_v = \{i,j,k\}$, if $p^*_{v,r} = 0$ then $p_{v,i}^* = p_{v,j}^* = p_{v,k}^* = 1$. Remember that each defender resource $r_j$ among the $m$ (non-red) defender resources can only protect a node of type $j$.  This means the $m$ non-red nodes in these $m/3$ sub-graphs are from $m$ different types. As a result, we obtain the satisfied cover for the original 3-Cover problem --- each subset belonging to this cover has three items that correspond to the three non-red nodes of one of the above $m/3$ sub-graphs.   
\end{proof}
\section{Gradient-based Binary Search Algorithm}
\label{sec:handle-exponential}
Overall, the problem of finding an optimal interdiction strategy for the defender as formulated in \eqref{prob:main-homo} is computationally challenging since the objective in \eqref{prob:main-homo} not only involves an exponential number of paths in the network but also is non-convex. We first introduce a new gradient-based binary search algorithm to solve \eqref{prob:main-homo} efficiently. Our algorithm has three key steps: (i) We use binary search to reduce the original fractional \eqref{prob:main-homo} to a simpler non-fractional problem; (ii) We present a non-trivial compact representation of the objective function based on the creation of a dynamic program, which handles an exponential number of paths involved in the original objective formulation; and (iii) We apply a gradient ascent-based method to efficiently solve the resulting compact optimization problem. Note that since our problem is non-convex, our gradient-based algorithm does not guarantee a global optimal solution. Therefore, we later introduce new guaranteed approximate solutions for \eqref{prob:main-homo} in Section~\ref{sec.approx.solution}. 

Essentially, we re-write the objective of \eqref{prob:main-homo} as follows:
\begin{align}
\cF^l(\bx) &= \sum_{s\in\cL}r^l(s, x_s)P^f_s(\bx) = \sum_{s\in\cL} \sum_{\substack{\tau \in\Omega\\ \tau \ni s}}r^l(s, x_s) P(\tau) = \frac{  \sum_{s\in\cL}\sum_{\tau \ni s}  r^l(s, x_s) \exp\left(\frac{U(\tau; \bx)}{\mu}\right)}{\sum_{\tau' \in \Omega} \exp\left(\frac{U(\tau';\bx)}{\mu}\right)}.
\label{eq:trajectory_formulation}
\end{align}
$\cF^l(\bx)$ has a fractional non-convex form. A typical way to simplify this structure is to use the Dinkelbach transform  and a binary search algorithm \citep{dinkelbach1967nonlinear} to convert the original problem into a sequence of simpler ones. 
We use binary search to write \eqref{prob:main-homo} equivalently as:
\begin{align}
&\max_{\delta}\left\{\delta\Big|\: \exists \bx \text{ s.t. } \cF^l(\bx) \geq \delta\right\} 
\Leftrightarrow \max_{\delta}\left\{\delta\Big|\: \max_{\bx} \left\{g(\bx,\delta)\right\} \geq 0\right\}, \nonumber
\end{align}
where the resulting objective:
\begin{align}
    g(\bx,\delta) =& \sum_{s\in\cL}\sum_{\substack{\tau\in \Omega \\ \tau \ni s}}  r^l(s, x_s) \exp\left(\frac{U(\tau;\bx)}{\mu}\right) - \delta \left(\sum_{\tau'\in\Omega} \exp\left(\frac{U(\tau';\bx)}{\mu}\right)\right)\label{eq:g-func}
\end{align}
Overall, $g(\bx,\delta)$ is still non-convex, but no longer fractional. Since $g(\bx,\delta)$ is differentiable, this maximization problem can be solved for a local maximum by a gradient-based method. One of the key challenges is the computation of $g(\bx,\delta)$, which, if done naively, would require enumerating exponentially many paths on $\Omega$. In the following, we show that $g(\bx,\delta)$ has a compact form, which allows us to compute $g(\bx,\delta)$ and its gradient efficiently via dynamic programming.   

\subsection{Compact Representation: Handling Exponential Numbers of Paths }
\label{sec:handling-exp-paths}
For a compact representation of $g(\bx,\delta)$, we introduce the following new terms:
\begin{align*}
    &Z_s = \sum_{\tau \in \Omega^{s_d}(s)} \exp\left(\frac{U(\tau;\bx)}{\mu}\right),\:\forall s\in \cS, &Y^s_{s'} = \sum_{\tau \in \Omega(s',s)} \exp\left(\frac{U(\tau;\bx)}{\mu}\right)
\end{align*}
where $\Omega^{s_d}(s)$ is the set of all {paths} from $s$ to the destination $s_d$ and $\Omega(s',s)$ is the set of all {paths} going from $s'$ to $s$, for any $s,s'\in \cS$. 

Given these new terms, the objective $g(\bx, \delta)$ can be re-formulated as follows: 
\begin{equation}
g(\bx, \delta) = \sum\nolimits_{s\in \cL}r^l(s, x_s) Y^s_{s_o}Z_s - \delta Z_{s_o}, \label{eq:gobj}
\end{equation}
where $s_o$ is the origin and $\Omega(s_o,s)$ is the set of all {paths} from $s_o$ to $s$. 

Although these new terms still involve exponentially many {paths} in $\Omega^{s_d}(s)$ and $\Omega(s',s)$ for all $s, s'$, we show that these terms can be computed efficiently via dynamic programming. Essentially, $Z_s$, $\forall s\in \cS$, can be computed recursively as follows:
\[
Z_s = \begin{cases}
\sum_{s'\in N(s)}\exp\left(\frac{v(s;\bx)}{\mu}\right) Z_{s'} & \text{ if }s \neq s_d \\
  1 & \text{ if } s = s_d,
\end{cases}
\]
where {$N(s)= \{s'\in \cS|
~ (s,s') \in \cA\}$}, is the set of possible next nodes that can be reached in one hop from node $s\in \cS$.  
Let $\bM$ be a matrix of size $|\cS\times \cS|$ with entries $M_{ss'} =  \exp\left(\frac{v(s|\bx)}{\mu}\right)$ for all $s\in \cS, s'\in N(s)$, then $\bZ = \{Z_s,~s\in \cS\}$ is a solution to the following linear system $\bZ = \bM\bZ +\bb$, where $\bb$ is a vector of size $|\cS|$ with zero entries except $b_{s_d} = 1$.  

Furthermore, $\bY^s = \{Y^s_{s'}\}_{s'}$ can also be defined recursively as follows:
\[
Y^s_{s'} = \begin{cases}
\sum_{s'' \in N(s')} \left(\exp\left(\frac{v(s';\bx)}{\mu}\right)\right) Y_{s''} & \text{ if }s' \neq s \\
  1 & \text{ if } s' = s.
\end{cases}
\]
Analogous to the technique above for $\bZ$, it can be seen that $\bY^s$
 is a solution to the linear system $\bY^s = \bM \bY^s + \bb^s$, where $\bb^s$ is of size $|\cS|$ with zeros everywhere except $b^s_s = 1$. 
 
Since $\bY^s$ and $ \bZ$ are solutions to the systems $\bY^s = \bM \bY^s + \bb^s$ and $\bZ = \bM\bZ +\bb$, respectively, $\forall s\in \cS$, the objective $g(\bx,\delta)$ can be computed via solving $|\cL|+1$ system of linear equations.
Finally, we see that all the above linear systems rely on the common matrix $\bM$. We can group them all into only one linear system.
Let $\bH$ be a matrix of size $(|\cS|)\times (|\cL|+1)$ in which the 1st to $|\cL|$-th columns are vectors $\bY^s$, $s\in \cL$ and the last column is $\bZ$. Let $\bB$ be a matrix of size $(|\cS|)\times (|\cL|+1)$ in which the 1st to $|\cL|$-th columns are vectors $\bb^s$, $s\in \cL$ and the last column is $\bb$. We see that $\bH$ is a solution to the linear system $(\bI-\bM)\bH = \bB$. Thus, in general, we can solve only one linear system to obtain all $\bY^s$ and $\bZ$. This way should be scalable when the size of $\cL$ increases. 
\subsection{Gradient-based Optimization}
As mentioned previously, we aim at employing the gradient-based approach to solve the following optimization problem (as the result of applying binary search to the original problem \eqref{prob:main-homo}):
\begin{align*}
    & \max_{\bx} \left\{g(\bx,\delta)\right\}
\end{align*}
In order to do so, the core is to compute the gradient $\left\{\frac{\partial g(\bx,\delta)}{\partial x_s}\right\}$. According to Equation~\ref{eq:gobj}, this gradient computation requires differentiating through the matrices $\bZ$ and $\{\bY^s\}$ (or equivalently, differentiating through the matrix $\bH$). We first present our Proposition~\ref{prop.2}:
 
\begin{proposition}\label{prop.2}
If the network is cycle-free, then $(\bI-\bM)$ is invertible.
\end{proposition}
Our Proposition~\ref{prop.2} shows that $(\bI-\bM)$ is invertible, allowing us to compute the matrix $\bH$: 
\[
\bH = (\bI-\bM)^{-1}\bB. 
\]
By taking the derivatives of both sides w.r.t $x_j$, $j\in\cL$, we obtain:
\begin{align}
\bJ^{\bH,j} &= (\bI-\bM)^{-1} \bJ^{\bM,j} (\bI-\bM)^{-1} \bB = (\bI-\bM)^{-1} \bJ^{\bM,j} \bH,\;\forall j\in \cL,\nonumber 
\end{align}
where $\bJ^{\bH,j}$ and $\bJ^{\bM,j}$ are the gradient matrices of $\bH$ and $\bM$ w.r.t $x_j$, i.e.,  $\bJ^{\bH,j}$ is a matrix of size $|\cS|\times (|\cL|+1)$ with entries $\bJ^{\bH,j}_{ss'} = \partial H_{ss'}/\partial x_j$, and  
$\bJ^{\bM,j}$ is a matrix of size $(|\cS|\times|\cS|)$ with entries $\bJ^{\bM,j}_{ss'} = \partial M_{ss'}/\partial x_j$, for any $j\in \cL$. Let $\bR^l(\bx)$ be a matrix of size $1 \times |\cL|$ with entries $r^l(s, x_s)$ for $s \in \cL$. We use $A_{S,T}$ to denotes a sub-matrix of $A$ which uses the rows in set $S$ and columns in set $T$. If $S$ or $T$ is a singleton, e.g., $S = \{s_o\}$ or $T= |\cL| + 1$, then we write it as $s_o$ or $|\cL|+1$. 

As a result, we now can compute the required gradient as follows:
\begin{align}\nonumber
    \frac{\partial g(\bx,\delta)}{\partial x_s} &= (\bR^l(\bx)\circ \bJ^{\bH,s}_{s_o,\cL} + \bJ^{R,s} \circ \bH_{s_o,\cL}) \times \bH_{\cL,|\cL|+1}  \\& + (\bR^l(\bx)\circ \bH_{s_o,\cL}) \times \bJ^{\bH,j}_{\cL,|\cL|+1}  -\delta \bJ^{\bH}_{s_o,|\cL|+1},   \forall s\in\cL \label{eq:grad}
\end{align}
where $\circ$ denotes Hadamard product.
We summarize  the main steps to optimize $g(\bx,\delta)$ in Alg. \ref{alg:first}.
\begin{algorithm}[htb]
\DontPrintSemicolon
\caption{\textit{Maximizing  $g(\bx, \delta)$}\label{alg:first}}  
\While{not converged}{
Given $\bx$, solve the system $\bH = (\bI-\bM)^{-1}\bB$ and $\bJ^{\bH,j} = (\bI-\bM)^{-1} \bJ^{\bM,j} \bH$ for all $j$ \;
Compute $g(\bx,\delta)$ and $\frac{\partial g(\bx,\delta)}{\partial x_s}$ using Eq.~\ref{eq:gobj},~\ref{eq:grad}. \;
Update $\bx$ using a projected gradient method
}
\end{algorithm}


\begin{remark}The above gradient-based approach only guarantees a local optimum. The complexity is determined by the matrix inversion or by solving a linear equation system, which, in worst case, is in $O(|\cS|^3)$. The gradient descent loop runs $O(1/\epsilon)$ to provide an additive $\epsilon$ approximation. Thus, the total complexity is $O((1/\epsilon) |\cS||\cA|)$. In practice, the gradients can be found using auto differentiation techniques, providing significantly more speed-up than directly using the formula for $\bJ^{\bH, j}$.
\end{remark}
\section{Guaranteed Approximate Solutions}
\label{sec.approx.solution}
As mentioned previously, while we can handle and exponential number of paths, our gradient-based method cannot guarantee to find a global optimal solution for \eqref{prob:main-homo} due to its non-convexity. In this section, we present results on top of the approach in Alg.~\ref{alg:first} that provides approximation guarantees. For this approximation, we need assumptions on the utility function, namely, that the utilities have a linear form: $v(s;\bx) =w^f_s x_s+ t^f_s $ and $r^l(s; \bx) = r^l(s, x_s) = w^l_s x_s + t^l_s$ for some constants $w^f_s, t^f_s, w^l_s, t^l_s$. We also assume that $w^f_s < 0$ and $w^l_s > 0$, i.e., 
more resources $x_s$ at node $s$ will lower adversary's utilities, and provide more  rewards to the defender. This setting is intuitive for security settings~\citep{yang2012computing,mai2022choices}. 

In the following, we first introduce a restricted interdiction problem that can be solved \emph{optimally} in a tractable manner using
our efficient gradient descent-based method. We then present our important theoretical results on its solution connection with the original problem. We further show that the restricted interdiction problem can be solved substantially more efficiently with a guaranteed solution bound via the binary search approach.
\paragraph{\textbf{Restricted Interdiction Problem.}}Let $\Delta(s)$ be the set of {paths} that cross $s$ and do not cross any other {node} in $\cL$.  We consider the following restricted interdiction problem:
\begin{align}
     \max_{\bx } &\;  \widetilde{\cF}(\bx)  = \frac{\sum_{s\in \cL}\sum_{\tau\in\Delta(s)} r^l(s, x_s) \exp\left(\frac{U(\tau;\bx)}{\mu}\right)}{\sum_{s\in \cL}\sum_{\tau\in\Delta(s)}  \exp\left(\frac{U(\tau;\bx)}{\mu}\right)} \label{prob:approx-main} \tag{\sf Approx-OPT}\\
    \mbox{subject to} & \; \sum\nolimits_{s\in \cL_k}x_s\leq M_k,~\forall k\in[K], \text{ and } x_s \in[L^x,U^x],\;\forall s\in \cL. \nonumber
\end{align}
Intuitively, in \eqref{prob:approx-main}, the adversary's path choices are restricted to a subspace of {paths} in the network which only cross a single critical {node} in $\cL$. We denote by $\cX = \{\bx \in \bbR^{|\cL|}|\;\sum_{s\in \cL_k}x_s\leq M_k,~\forall k\in [K],~x_s \in[L^x,U^x],\;\forall s\in \cL \}$, the feasible set of the defender's interdiction strategies $\bx$. 
\subsection{Solution Relation with Original Problem \eqref{prob:main-homo}}
Given the formulation of restricted problem \eqref{prob:approx-main}, we now theoretically analyze its solution relation with our original problem \eqref{prob:main-homo} of finding an optimal interdiction strategy for the defender.
\subsubsection{\textbf{Main results}}
Given a critical node $s\in \cL$, let $\Delta^+(s)$ be the set of {paths} that cross $s$ and at least another critical {node} in $\cL$.  
Let $\beta_1$, $\beta_2>0$ such that:
\begin{equation}\label{eq:beta}
\begin{aligned}
    \beta_1 &= \max_{\bx}\max_{s\in \cL} \left\{\frac{\sum_{\tau \in \Delta^+(s)}  \exp\left(\frac{U(\tau;\bx)}{\mu}\right) }{\sum_{\tau \in \Delta(s)}  \exp\left(\frac{U(\tau;\bx)}{\mu}\right)}\right\}\\
        \beta_2 &= \max_{\bx} \left\{\frac{\sum_{\tau \in \bigcup_s\{\Delta^+(s)\}}  \exp\left(\frac{U(\tau;\bx)}{\mu}\right) }{\sum_{\tau \in \bigcup_s\{\Delta(s)\}}  \exp\left(\frac{U(\tau;\bx)}{\mu}\right)}\right\},
\end{aligned}
\end{equation}
Intuitively, $\beta_1$ and $\beta_2$ capture the maximum change in the adversary's behavior when the adversary's path choice is limited to the restricted path space as in \eqref{prob:approx-main}.  The values of $\beta_1$ and $\beta_2$ are expected to be small if the cost of traveling between any two critical nodes in $\cL$ is large. That is, $\beta_1,\beta_2\rightarrow 0$ if $\sum_{\tau\in \Omega(s,s')} \exp\left(\frac{U(\tau;\bx)}{\mu}\right)\rightarrow 0$,
where $\Omega(s,s')$ is the set of all {paths} going from $s$ to $s'$, for any $s,s'\in \cL$. Our main theoretical results are presented in Theorems ~\ref{thm.restricted.relation.1} and \ref{thm.restricted.relation.2}.

\begin{theorem}\label{thm.restricted.relation.1}
Let $\bx^*$ be optimal to \eqref{prob:approx-main}: $\max_{\bx\in \cX} \widetilde{\cF}(\bx)$ and $\kappa = \max_{\bx\in\cX}\sum_{s\in \cL}|r^l(s,x_s)|$ the maximal absolute reward the defender can possibly archive at a critical node, we have:
\[
\cF^l(\bx^*) \geq \frac{\max_{\bx \in \cX}\{ \cF^l(\bx)\} }{(1+\beta_1)(1+\beta_2)} - \kappa\frac{\beta_1+\beta_2+\beta_1\beta_2}{(1+\beta_1)(1+\beta_2)}.
\]
\end{theorem}
Note that $\max_{\bx \in \cX}\{ \cF^l(\bx)\}$ is the original \eqref{prob:main-homo} to find an optimal interdiction strategy. As stated previously, when the cost of traveling between any two critical {nodes} is high, $\beta_1$ and $\beta_2$ get closer to zero, meaning the RHS of the inequality in Theorem~\ref{thm.restricted.relation.1} will reach the optimal solution of \eqref{prob:main-homo}. 

Now, if we can only solve \eqref{prob:approx-main} approximately, we can still estimate the approximation bound for the original \eqref{prob:main-homo}, as shown in Theorem~\ref{thm.restricted.relation.2}.
\begin{theorem}\label{thm.restricted.relation.2}
Given $\epsilon>0$, let $\bx^*$ be a solution such that  $\widetilde{\cF}(\bx^*) \geq (1-\epsilon)\max_{\bx} \widetilde{\cF}(\bx)$, we have:
\[
{\cF^l(\bx^*)} \geq  \frac{(1-\epsilon) {\max_{\bx} \{\cF^l(\bx)\}}}{(1+\beta_1)(1+\beta_2)} - \kappa\frac{\epsilon + \beta_1+\beta_2+\beta_1\beta_2}{(1+\beta_1)(1+\beta_2)}.
\]
Moreover, if $\bx^*$ is an approximate solution with an additive error $\epsilon>0$, i.e.,  $\widetilde{\cF}(\bx^*) \geq \max_{\bx} \widetilde{\cF}(\bx) -\epsilon$, we obtain the following bound:
\[
\cF^l(\bx^*) \geq \frac{1}{\eta} \max_{\bx} \{\cF^l(\bx)\} -\kappa \frac{\eta-1}{\eta}, 
\]
where  $\eta = (1+\beta_1)(1+\beta_2)\left(1+\frac{\epsilon}{\kappa + \min_{\bx\in\cX} \widetilde{\cF}(\bx)}\right)$.
\end{theorem}
\subsubsection{\textbf{Proof Sketches. }} The proof of these theorems are based on two important lemmas, as explained below. Lemma~\ref{lemma.1} only applies when all the defender's rewards $r^l(s,x_s)$ are non-negative. We then handle the  general case in Lemma~\ref{lm:choose-kappa}. Intuitively, these two lemmas show relations in terms of the defender's utilities (aka. objective functions $F^l(\bx)$ and $\widetilde{F}(\bx)$) between the original problem \eqref{prob:main-homo} and the restricted problem \eqref{prob:approx-main} for any given defender's interdiction strategy $\bx$.
\begin{lemma}\label{lemma.1}
If $r^l(s,x_s) \geq 0$ for any $\bx\in\cX$, then for any $\bx$, 
\[
\frac{1}{1+\beta_2} \widetilde{\cF}(\bx)\leq \cF^l(\bx)\leq (1+\beta_1) \widetilde{F}(\bx).
\]
\end{lemma}
The case that $r^l(s,x_s)$ would take negative values is more challenging to handle. The following lemma gives general inequalities for such a situation. 
\begin{lemma}
\label{lm:choose-kappa}
If we choose 
$
\kappa =\sum_{s\in \cL}\max_{\bx}\left|r^l(s,x_s) \right|,
$
then for any $\bx\in\cX$, 
\[
\frac{1}{1+\beta_2} \left(\widetilde{\cF}(\bx)+\kappa\right)\leq \cF^l(\bx) +\kappa \leq (1+\beta_1) \left(\widetilde{\cF}(\bx)+\kappa\right).
\]
\end{lemma}
The proofs of Theorem~\ref{thm.restricted.relation.1} and Theorem~\ref{thm.restricted.relation.2} now follow readily from the above lemmas.

\begin{proof}[\text{Completing Proof of Theorem~\ref{thm.restricted.relation.1}}]
According to Lemma \ref{lm:choose-kappa}, we obtain:
\begin{align}\label{eq:lm36-eq1}
\frac{1}{1+\beta_2}\max_{\bx} \left(\widetilde{\cF}(\bx) +\kappa \right) &\leq \max_{\bx} \left(\cF^l(\bx) +\kappa  \right)\leq (1+\beta_1) \max_{\bx} \left(\widetilde{\cF}(\bx)+\kappa\right),  
\end{align}
leading to the following chain of inequalities:
\begin{align*}
    \frac{1}{1+\beta_2}\left(\widetilde{\cF}(\bx^*)+\kappa\right)&\leq  \left(\cF^l(\bx^*)+\kappa\right) \leq  \max_{\bx} (\cF^l(\bx)+\kappa) \leq (1+\beta_1) \left(\widetilde{\cF}(\bx^*)+\kappa\right). 
\end{align*}
where $\bx^*$ is the optimal defense strategy solution to \eqref{prob:approx-main}. Since $\widetilde{\cF}(\bx^*)+\kappa, \cF^l(\bx^*)+\kappa$, and $\max_{\bx} (\cF^l(\bx)+\kappa)$ are all positive, we have:
\begin{align*}
    \frac{\cF^l(\bx^*)+\kappa}{\max_{\bx} (\cF^l(\bx)+\kappa)} &\geq \frac{ \frac{1}{1+\beta_2}\left(\widetilde{\cF}(\bx^*)+\kappa\right)}{(1+\beta_1) \left(\widetilde{\cF}(\bx^*)+\kappa\right)}= \frac{1}{(1+\beta_1)(1+\beta_2)}\\
    \text{which implies: } 
     \cF^l(\bx^*) +\kappa &\geq \frac{\max_{\bx \in \cX}\{ \cF^l(\bx)\}}{(1+\beta_1)(1+\beta_2)}  + \frac{\kappa}{(1+\beta_1)(1+\beta_2)}\\
    \text{or equivalently, } 
    \cF^l(\bx^*)& \geq \frac{\max_{\bx \in \cX}\{ \cF^l(\bx)\}}{(1+\beta_1)(1+\beta_2)}  - \kappa\frac{\beta_1+\beta_2+\beta_1\beta_2}{(1+\beta_1)(1+\beta_2)}
\end{align*}
which concludes our proof.
\end{proof}
\begin{proof}[\text{Completing Proof of Theorem~\ref{thm.restricted.relation.2}}]
If $\bx^*$ is the interdiction strategy solution such that $\widetilde{\cF}(\bx^*) \geq (1-\epsilon)\max_{\bx} \widetilde{\cF}(\bx)$, then by using \eqref{eq:lm36-eq1}, we obtain the following chain of inequalities:
\begin{align}
\frac{1}{1+\beta_2} \left(\widetilde{\cF}(\bx^*) +\kappa \right) &\leq \frac{1}{1+\beta_2}\max_{\bx} \left(\widetilde{\cF}(\bx) +\kappa \right) \leq \max_{\bx} \left(\cF^l(\bx) +\kappa  \right)\\\nonumber
&\leq (1+\beta_1) \max_{\bx} \left(\widetilde{\cF}(\bx)+\kappa\right)\leq 
(1+\beta_1)  \left(\frac{1}{1-\epsilon}\widetilde{\cF}(\bx^*)+\kappa\right) \\\nonumber
&\leq \frac{1+\beta_1}{1-\epsilon}(\widetilde{\cF}(\bx^*)+\kappa)
\end{align}
Thus, 
\begin{align}\nonumber
    &\frac{1}{1+\beta_2} \left(\widetilde{\cF}(\bx^*) +\kappa \right) \leq {\cF^l}(\bx^*) +\kappa  \leq \max_{\bx} \left(\cF^l(\bx) +\kappa  \right) \leq \frac{1+\beta_1}{1-\epsilon}(\widetilde{\cF}(\bx^*)+\kappa)\\
    \implies &\frac{\cF^l(\bx^*)+\kappa}{\max_{\bx} \{\cF^l(\bx)+\kappa\}} \geq  \frac{1-\epsilon}{(1+\beta_1)(1+\beta_2)},
\end{align}
which further leads to:
\[
{\cF^l(\bx^*)} \geq  \frac{(1-\epsilon) {\max_{\bx} \{\cF^l(\bx)\}}}{(1+\beta_1)(1+\beta_2)} - \kappa\frac{\epsilon + \beta_1+\beta_2+\beta_1\beta_2}{(1+\beta_1)(1+\beta_2)},
\] 
 as desired.
 
For the case of additive error  $\widetilde{\cF}(\bx^*) \geq \max_{\bx} \widetilde{\cF}(\bx)-\epsilon$, similarly, we can write:
\begin{align}
\frac{1}{1+\beta_2} \left(\widetilde{\cF}(\bx^*) +\kappa \right) &\le
{\cF^l}(\bx^*) +\kappa 
 \leq\max_{\bx} \left(\cF^l(\bx) +\kappa  \right)
\nonumber\\
&\leq (1+\beta_1) \max_{\bx} \left(\widetilde{\cF}(\bx)+\kappa\right) \leq 
(1+\beta_1)  \left(\widetilde{\cF}(\bx^*) +\epsilon+\kappa\right), 
\nonumber 
\end{align}
which yields:
\begin{align*}
   \frac{\cF^l(\bx^*)+\kappa}{\max_{\bx} \{\cF^l(\bx)+\kappa\}} &\geq  \frac{1}{(1+\beta_1)(1+\beta_2)} \frac{1}{1+\epsilon/(\widetilde{\cF}(\bx^*) +\kappa )}\geq \frac{1}{\eta}\\
   \Rightarrow  \cF^l(\bx^*) &\geq \frac{1}{\eta} \max_{\bx} \{\cF^l(\bx)\} -\kappa \frac{\eta-1}{\eta}, 
\end{align*}
as desired, which concludes our proof.
\end{proof}

\subsection{Solving the Restricted Interdiction Problem  \eqref{prob:approx-main}}
In order to solve the restricted problem, we also propose to apply the binary search approach. The resulting sub-problem at each binary search step can be formualted as follows:
\begin{align}
\nonumber &\max_{\bx \in\cX}\; \widetilde{g}(\bx,\delta) = \sum_{s\in\cL}\sum_{\tau \in \Delta(s)}  r^l(s, x_s) \exp\left(\frac{U(\tau;\bx)}{\mu}\right)  - \delta \left(\sum_{s\in\cL}\sum_{\tau \in \Delta(s)} \exp\left(\frac{U(\tau;\bx)}{\mu}\right)\right).\label{eq:app-sub}\tag{\sf sub-Approx}  
\end{align}
 
The rest of this section will devote to presenting our theoretical results on the key underlying property of \eqref{eq:app-sub}, as well as new exact/approximate solutions for solving \eqref{eq:app-sub}.

\subsubsection{Unimodality of \eqref{eq:app-sub}}
In the following, we present Theorem~\ref{prop:c1} which shows that we can use a gradient-based method to obtain the unique global optimal solution to \eqref{eq:app-sub}.
\begin{theorem}
\label{prop:c1}
The \eqref{eq:app-sub} problem is unimodal; that is, any local optimal solution $\bx^*$ of \eqref{eq:app-sub} is also globally optimal.
\end{theorem}
\begin{proof}[Proof Sketch]
The proof can be divided into two major steps: 

\paragraph{\textbf{Step 1: showing that \eqref{eq:app-sub} can be converted into a (strictly) convex optimization problem.}}

This step is done via variable transformation. Recall that the adversary utility $v(s;\bx) =v(s, x_s)=w^f x_s+ t^f $ and defender utility $r^l(s;\bx) = r^l(s,x_s) = w^l x_s + t^l$ where $w^f < 0$ and $w_l > 0$. Essentially, we introduce a new variable $y_s = \exp\left(\frac{v(s, x_s)}{\mu}\right)$ for all critical nodes $s\in\cL$. We will show that the objective $\widetilde{g}(\bx,\delta)$ of \eqref{eq:app-sub} can be rewritten as a strictly concave function of $\{y_s\}$.

For each trajectory $\tau \in \Delta(s)$, let $V^s(\tau) = \sum\limits_{s'\in\tau,\; s'\neq s} v(s'; \bx)$, or equivalently, $V^s(\tau) = U(\tau;\bx) - v(s;\bx)$, which is the accumulated adversary utility over every node on $\tau$ except  node $s$. We see that $V^s(\tau)$ is independent of the defender coverage probability $x_s$ at every critical node $s\in\cL$ (by definition of $\Delta(s)$). Therefore, we can reformulate the objective $\widetilde{g}(\bx,\delta)$ of \eqref{eq:app-sub} as follows:
\begin{align*}
    \widetilde{g}(\bx,\delta) = \sum_{s\in\cL} r^l(s, x_s) \exp\left(\frac{v(s;\bx)}{\mu}\right) H(s) - \delta \left(\sum_{s\in\cL}\exp\left(\frac{v(s;\bx)}{\mu}\right) H(s)\right).
\end{align*}
where $H(s) =\sum_{\tau \in \Delta(s)} \exp(V^s(\tau)/\mu)$.   
We thus can write  $\widetilde{g}(\bx,\delta)$ as function of $\by$ as follows: 
\begin{align}
\widetilde{g}(\by,\delta)&= \sum_{s\in\cL}\left(\left(\mu \log (y_s)-t^f\right) \frac{w^l}{w^f} + t^l\right) y_s H(s) -\delta\left(\sum_{s\in\cL} y_s H(s)\right) \nonumber \\
&=\sum_{s\in\cL} \mu\frac{w^l}{w^f} H(s)  \log(y_s)y_s - y_sH(s)\left(\mu\frac{t^fw^l}{w^f} +t^l+\delta\right).
\end{align}
Since  $w^f/w^l\leq  0$, it can be shown that each component $\mu\frac{w^l}{w^f} H(s)  \log(y_s)y_s$ is concave in $y_s$, thus
$\widetilde{g}(\by,\delta)$ is strictly concave in $\by$ for all critical nodes $s$. Moreover, {for any $k\in[K]$}, the constraint $\sum_{s\in \cL_k} x_s \leq  M_k$ becomes $\sum_{s\in \cL_k}  \frac{\mu \log (y_s)-t^f}{w^f} \leq M_k$, which is convex since $w^f <0$.

\paragraph{\textbf{Step 2: proving global optimality via the KKT condition correspondence with variable transformation}} Under the variable transformation as presented in \textbf{Step 1}, for notational convenience, let us define $\widehat{x}_s(\cdot):\bbR\rightarrow\bbR$ and $\widehat{y}_s(\cdot):\bbR\rightarrow\bbR$ such that: 
\begin{align*}
    &\widehat{y}_s(x_s) = \exp\left(\frac{v(s, x_s)}{\mu}\right)\\
    &\widehat{x}_s(y_s) = \frac{\mu \log (y_s)-t^f}{w^f},
\end{align*}
i.e., the mappings from $\bx_s$ to $\by_s$ and vice-versa. 

Recall that the feasible strategy space of the defender $\cX=\left\{\bx: \sum_{s \in \cL_k} x_s \leq M_k,  \forall k, x_s\in [L^x,U^x]\right\}$. We thus can write the Lagrange dual of \eqref{eq:app-sub} as follows:
\begin{align*}
    L^g(\bx,\bgamma,\bbeta^1,\bbeta^2)& =  \widetilde{g}(\bx, \delta) - \sum_k \gamma_k \left(\sum_{s \in \cL_k} x_s - M_k \right) - \sum_{s}\eta^1_s (x_s - U^x) + \sum_{s}\eta^2 (x_s-L^x).
\end{align*}
Since $\bx^*$ is a local optimal solution for \eqref{eq:app-sub}, the KKT conditions imply that there are dual $\bgamma^*,\bbeta^{1*}, \bbeta^{1*}\geq 0$ such that the following constraints are satisfied:
\begin{equation}\label{eq:KKT-x}
    \begin{cases}
    \frac{\widetilde{g}(\bx^*, \delta)}{\partial x_s} - \gamma^*_k  - \eta^{1*}_s + \eta^{2*}_s = 0, \mbox{ where $k$ such that } s \in \cL_k\\
    \gamma^*_k \left(\sum_{s \in \cL_k} x_s - M_k\right) = 0, \; \forall k\\
    \eta^{1*}_s (x^*_s - U^x) = 0\\
    \eta^{2*}_s (x^*_s - L^x) = 0\\
    L^x\leq x^*_s \leq U^x\\
    \sum_{s \in \cL_k} x_s - M_k, \; \forall k
    \end{cases}
\end{equation}
By the variance transformation $y_s   = \exp\left(\frac{(w^fx_s +t^f)}{\mu}\right)$, let $y^*_s = \exp\left(\frac{(w^fx^*_s +t^f)}{\mu}\right)$
and $x_s = \widehat{x}_s(y_s) = \frac{\mu \log (y_s)-t^f}{w^f}$ 
for all $s\in \cL$, we can write  \eqref{eq:KKT-x} equivalently  as:
\begin{equation}\label{eq:KKT-y}
    \begin{cases}
    \frac{\widetilde{g}(\bx^*, \delta)}{\partial x_s} \frac{\partial \widehat{x}_s(y^*_s)}{\partial y_s} - \gamma^*_k \frac{\partial \widehat{x}_s(y^*_s)}{\partial y_s}  - \eta^{1*}_s \frac{\partial \widehat{x}_s(y^*_s)}{\partial y_s} + \eta^{2*}_s \frac{\partial \widehat{x}_s(y^*_s)}{\partial y_s} = 0\\
    \gamma^*_k \left(\sum_{s\in \cL}  \frac{\mu \log (y^*_s)-t^f}{w^f} -  M\right) = 0, \; \forall k\\
    \eta^{1*}_s \left(\frac{\mu \log (y^*_s)-t^f}{w^f}  - U^x)\right) = 0\\
    \eta^{2*}_s \left( \frac{\mu \log (y^*_s)-t^f}{w^f}  - L^x \right) = 0\\
    L^x\leq  \frac{\mu \log (y^*_s)-t^f}{w^f}  \leq U^x\\
    \sum_{s\in \cL_k}  \frac{\mu \log (y^*_s)-t^f}{w^f} \leq M_k, \; \forall k.
    \end{cases}
\end{equation}
The first condition of \eqref{eq:KKT-y} can be written equivalently as follows: 
\[
\frac{\widetilde{g}(\by^*, \delta)}{\partial y_s} - \gamma^*_k \bbI[s\in \cL_k] \frac{\partial \widehat{x}_s(y^*_s)}{\partial y_s}  - \eta^{1*}_s \frac{\partial \widehat{x}_s(y^*_s)}{\partial y_s} + \eta^{2*}_s \frac{\partial \widehat{x}_s(y^*_s)}{\partial y_s} = 0,
\]
where $\bbI[\cdot]$ is the indicator function. This
implies that $\by^*,\bgamma^*,\bbeta^{1*}, \bbeta^{1*}$ also satisfy the KKT conditions of the following (strictly) convex optimization problem (i.e., the resulting problem of variable transformation discussed in \textbf{Step 1}).
\begin{align}
     \max_{\by } &\quad\quad  \widetilde{g}(\by,\delta) \label{prob:convex-1}\\
    \mbox{subject to} & \quad\quad \sum_{s\in \cL_k}\widehat{x}_s(y_s)\leq M_k, \; \forall k \text{ and } \widehat{x}_s(y_s) \in [L^x,U^x].\nonumber
\end{align}
Thus, $\by^*$ is the unique global optimal solution to \eqref{prob:convex-1}, which also means that $\bx^*$ is also the global optimal solution to \eqref{eq:app-sub} as desired. 
\end{proof}

\subsubsection{Exact Solution} 
Based on Theorem~\ref{prop:c1}, we can use a gradient-based method to obtain the unique global optimal solution to \eqref{eq:app-sub}. The computational challenge is that the objective $\widetilde{g}(\bx,\delta)$ still involves exponentially many {paths}. Our idea is to decompose $\widetilde{g}(\bx,\delta)$ into multiple terms (each term corresponds to a critical node $s\in\cL$) --- which can be computed using dynamic programming. Essentially, we create new graphs by keeping a node $s\in \cL$ and remove every other {nodes} in $\cL$. By doing this, we can exactly handle $\Delta(s)$ using our approach in Section \ref{sec:handling-exp-paths}. To facilitate our exposition,  let $\cG(s)$ be the sub-graph created from the graph $\cG$ by deleting all nodes in the set $\cL$ except node $s$. 
Since we will be dealing with several graphs, henceforth we denote all {paths} in any arbitrary graph $\cG$ as $\Omega(\cG)$. The proposition below shows that $\widetilde{g}(\bx,\delta)$ can be decomposed into terms that can be efficiently computed based on sub-graphs $\cG(s)$, $\forall s\in\cS$.
\begin{proposition}
$\widetilde{g}(\bx,\delta)$ can be written as follows:
\begin{align}\small
\label{eq:exact-tg} \sum_{s\in\cL} \left(\sum_{\substack{\tau \in \Omega(\cG(s))\tau \ni s}}\left[  r^l(s, x_s) \exp\left(\frac{U(\tau;\bx)}{\mu}\right) - \delta \exp\left(\frac{U(\tau;\bx)}{\mu}\right)\right]\!\right) 
\end{align}
\end{proposition}
Given any sub-graph $\cG(s)$, the terms in $\widetilde{g}(\bx,\delta)$ (as well as their gradients)
can be computed by solving a system of linear  equations, similarly as the approach described in in Algorithm~\ref{alg:first} of Section~\ref{sec:handling-exp-paths}. We solve the overall problem in the Algorithm~\ref{algo:main-exact} below. 


\begin{algorithm}[htp]
\DontPrintSemicolon
\caption{\textit{Solving \eqref{prob:main-homo} through the restricted interdiction problem \eqref{prob:approx-main}}\label{algo:main-exact}} 
- Create $\cG(s)$ for all $s\in \cL$ \;
- Using binary search to solve \eqref{prob:approx-main}. At each step of the binary search, optimize $\max_{\bx \in \cX} \widetilde{g}(\bx, \delta)$ using a gradient-based method utilizing the reformulation in Equation~\eqref{eq:exact-tg}.\;
- Re-solve \eqref{prob:main-homo} with the original network to improve the
solution given by the above step.
\end{algorithm}

\subsubsection{{Efficient Approximate Solution}} 
To recap, in previous sections, we provide: (i) a new algorithm to solve the original interdiction problem \eqref{prob:main-homo} --- this algorithm is efficient but only guarantees a local optimal solution; and (ii) a new algorithm (Algorithm~\ref{algo:main-exact}) to solve the restricted interdiction problem \eqref{prob:approx-main} --- this algorithm provides a global optimal solution for \eqref{prob:approx-main}, and as a result, guarantees a bounded approximate solution for \eqref{prob:main-homo} (Theorems~\ref{thm.restricted.relation.1}\&\ref{thm.restricted.relation.2}). However, when $|\cL|$ is large, this algorithm might not scale in a practical implementation as it needs $|\cL|$ sub-graphs. 

In this section, we thus propose a new algorithm which is both efficient and guarantees a bounded approximate solution for \eqref{prob:main-homo}. Our main ideas can be summarized as follows: (a) we identify a graph modification and solve \eqref{prob:main-homo} with the modified graph using the algorithm described in Section~\ref{sec:handle-exponential}; (b) despite that (a) does not guarantee a global optimal solution for \eqref{prob:main-homo}, we theoretically shows that this  resulting solution is a bounded approximate solution for \eqref{prob:approx-main} (Theorem~\ref{prop.temp}); and (c) finally, by leveraging findings in Theorem~\ref{thm.restricted.relation.2} and Theorem~\ref{prop.temp}, we obtain Theorem~\ref{th:overal-bound} showing that this resulting solution (obtained from (a)) is also a bounded approximate solution for \eqref{prob:main-homo}. We elaborate our ideas in the following.  

\paragraph{\textbf{Network modification.}}Essentially, we modify the original network $\cG$ by raising the costs of travelling between any pair of nodes in $\cL$ in such a way that $\beta_1$ and $\beta_2$ become arbitrarily small. We remark that, given any $\epsilon' >0$, we can always modify the travelling costs between pairs of nodes in $\cL$ to obtain a modified network $\cG'$ such that the following conditions holds:
\begin{align}
\max_{\bx}\max_{s\in \cL} \left\{\frac{\sum_{\tau \in \Delta^+(s)}  \exp\left(\frac{U(\tau;\bx)}{\mu}\right) }{\sum_{\tau \in \Omega, \tau \ni s}  \exp\left(\frac{U(\tau;\bx)}{\mu}\right)}\right\} &\leq \epsilon' \label{eq:eps-1}\\
        \max_{\bx}\left\{\frac{\sum_{\tau \in \bigcup_s\{\Delta^+(s)\}}  \exp\left(\frac{U(\tau;\bx)}{\mu}\right) }{\sum_{\tau \in\Omega}  \exp\left(\frac{U(\tau;\bx)}{\mu}\right)}\right\} &\leq \epsilon'. \label{eq:eps-2}
\end{align}
We remind that $\Delta^+(s)$ be the set of paths that cross $s$ and at least another node in $\cL$ and $\Omega$ is the set of all paths. We denote the objective of \eqref{eq:g-func} of the binary search step for \eqref{prob:main-homo} with respect to the modified network $\cG'$ as ${g}(\bx,\delta|\cG')$. We can optimize ${g}(\bx,\delta|\cG')$ by running gradient descent. The problem $\max_{\bx}{g}(\bx,\delta|\cG')$ is not convex, yet we can provide the following strong guarantee.

\paragraph{\textbf{Solution theoretical bounds.}}Let us first define:
\begin{align}\nonumber
    &\rho^s = \sum_{\substack{\tau\in \Omega,\tau\ni s}}   \exp\left(\frac{U(\tau;L^x \bbe)}{\mu}\right),\;\forall s\in \cL\\\nonumber
    & \rho = \sum_{\substack{\tau\in \Omega}}   \exp\left(\frac{U(\tau;L^x \bbe)}{\mu}\right)
\end{align}
where $\bbe$ is an all-one vector of size $|\cL|$. We remind that $L^x$ and $U^x$ are the lower and upper bounds on the resource coverage $x_s$ of the defender at every critical node $s\in \cL$. In addition, the utility functions for the adversary and defender are in the linear forms: $v(s; \bx) = v(s, x_s) = w^f_s x_s+ t^f_s $ and $r^l(s; \bx) = r^l(s, x_s) = w^l_s x_s + t^l_s$ for some constants $w^f_s, t^f_s, w^l_s, t^l_s$, where $w^f_s < 0$ and $w^l_s > 0$. 

\begin{theorem}\label{prop.temp}
If we run binary search to solve \eqref{prob:main-homo} with the modified network $\cG'$ and obtain $(\overline{\bx},\overline{\delta})$, the following performance bound is guarantee for the restricted problem \eqref{prob:approx-main}:
\begin{align}\label{lm:bound-approx}
    &\widetilde{\cF}(\overline{\bx}) \geq  \max_{\bx} \{\widetilde{F}(\bx)\} - \frac{\epsilon'( \cH+2\cC)}{\sum_{s\in\cL}\sum_{\tau \in \Delta(s)} \exp\left(\frac{U(\tau;U^x\bbe)}{\mu}\right)} 
\end{align}
where the constant $\cC = \left(\kappa \max_s\{\rho^s\} + \overline{\delta} {\rho}\right)$.
\end{theorem}
\begin{proof}[Proof Sketch]
Our proof relies on the two lemmas below. Lemma~\ref{lemma.temp} shows the relationship between the objective of binary search step $g(\bx,\delta|\cG')$ of \eqref{prob:main-homo} with modified $\cG'$ and the objective of binary search step $\widetilde{g}(\bx,\delta)$ of the restricted problem  \eqref{prob:approx-main} with original graph $\cG$.
\begin{lemma}\label{lemma.temp}
For any $(\bx, \delta)$, then we have:
\[
|g(\bx,\delta|\cG') - \widetilde{g}(\bx,\delta)|  \leq \epsilon'\left(\kappa \max_s\{\rho^s\} + \delta \rho\right). 
\]
\end{lemma}
Given $\delta$, Lemma~\ref{prop:gap-g-tilde} below shows that any local optimal solution to  $\max_{\bx}g(\bx,\delta|\cG')$ (i.e., the binary step of \eqref{prob:main-homo} with modified $\cG'$) is {in an $\cO(\epsilon')$ neighborhood of}  optimal solutions to \eqref{eq:app-sub}: $\max_{\bx}\widetilde{g}(\bx,\delta)$, i.e., the binary step of the restricted problem  \eqref{prob:approx-main} with original graph $\cG$. To prove the lemma, we show that, given any  local optimal solution $\overline{\bx}$ to $\max_{\bx}g(\bx,\delta|\cG')$, there is an unimodal function (i.e., becomes strictly concave after the change of variables as in the proof of Theorem \ref{prop:c1}) which is in an $\cO(\epsilon')$ neighborhood of $\widetilde{g}(\bx,\delta)$ and has $\overline{\bx}$  as a global optimal solution. 
\begin{lemma}
\label{prop:gap-g-tilde}
Let $\overline{\bx}$ be a local optimal  solution of  $\max_{\bx}g(\bx,\delta|\cG')$ for a given $\delta$, then we have:
\[
\max_{\bx}\{\widetilde{g}(\bx,\delta)\} - \widetilde{g}(\overline{\bx},\delta) \leq  
\epsilon' \cH, 
\]
\text{ where } $\cH = 2\left( \mu\rho^j + \frac{\kappa}{|w^f_j|} \max_{s}\{\rho^s\} +{\delta}   \rho\right)\exp\left(\frac{w^f_j(L^x-U^x)}{\mu}\right)$.
\end{lemma}

Let $(\bx^*,\delta^*)$ be the output of the bisection (aka binary search) to solve the restricted interdiction problem \eqref{prob:approx-main}. Based on the above two lemmas, we now try to bound the gap $|\overline{\delta}-\delta^*|$ (based on which we can bound $\max_{\bx\in \cX}\widetilde{\cF}(\bx) - \widetilde{\cF}(\overline{\bx})$ as we will explain later). 

First, since $(\overline{\bx},\overline{\delta})$ is an output of binary search for solving \eqref{prob:main-homo} with the modified network $\cG'$, we have $|g(\overline{\bx},\overline{\delta}|\cG')|\leq \xi$, where $\xi$ is a positive constant depending on the precision of the binary search. In fact, this constant can be made arbitrarily small and therefore, we can remove it from the rest of our proof for the sake of presentation --- that is, in the following, we simply consider $|g(\overline{\bx},\overline{\delta}|\cG')|=0$. Now according to Lemma 3, we have:
\begin{align}\label{noname.1}
    &|g(\overline{\bx},\overline{\delta}|\cG') - \widetilde{g}(\overline{\bx},\overline{\delta})|  \leq \epsilon'\left(\kappa \max_s\{\rho^s\} + \overline{\delta} \rho\right) =\epsilon' \cC
\end{align}

Since $|g(\overline{\bx},\overline{\delta}|\cG')|=0$ we have $|\widetilde{g}(\overline{\bx},\overline{\delta})|\leq \epsilon' \cC$. 

On the other hand, since  $(\bx^*,\delta^*)$ is the result of binary search for \eqref{prob:approx-main}, we have: $\widetilde{g}(\bx^*,\delta^*) = 0$ and $\delta^* = \max_{\bx\in \cX} \widetilde{F}(\bx)$.  
We denote by: 
$$\widetilde{K}(\delta) = \max_{\bx\in\cX}\widetilde{g}(\bx, \delta)$$
which is monotonic decreasing in $\delta$. In addition, $\widetilde{K}({\delta^*}) = \widetilde{g}(\bx^*,\delta^*) = 0$. 
According to Lemma \ref{prop:gap-g-tilde}:
\begin{equation}
\label{eq:th8-eq1}
 \widetilde{K}(\overline{\delta}) - \epsilon' 
\cH \le \widetilde{g}(\overline{\bx},\overline{\delta}) \le \widetilde{K}(\overline{\delta}). 
\end{equation}
As a result, we obtain the following chain of inequalities:
\begin{align}
    &\epsilon' \cH \geq |\widetilde{K}(\overline{\delta}) - \widetilde{g}(\overline{\bx},\overline{\delta})|\stackrel{(a)}{\geq} |\widetilde{K}(\overline{\delta})| - |\widetilde{g}(\overline{\bx},\overline{\delta})|\stackrel{}{\geq}|\widetilde{K}(\overline{\delta})| -\epsilon'\cC\label{eq:th8:eq3}\\
    \implies& |\widetilde{K}(\overline{\delta})| \leq \epsilon' (\cH+\cC)
\end{align}
where $(a)$ is due to the triangle inequality. We further have $\widetilde{K}({\delta^*}) = 0$, leading  to:  
$$|\widetilde{K}(\overline{\delta}) - \widetilde{K}(\delta^*)|\leq \epsilon' (\cH+\cC)$$
Since $\widetilde{K}(\delta)$ is differentiable in $\delta$,  the mean value theorem tells us that there is $\alpha\in  [\overline{\delta},\delta^*]$ 
such that:
\begin{equation}
\label{eq:lemma59-eq1}
|\widetilde{K}(\overline{\delta}) - \widetilde{K}(\delta^*)| = \widetilde{K}'(\alpha)|\overline{\delta}-\delta^*|\leq \epsilon' (\cH+\cC).
\end{equation}
We denote by $\widetilde{\bx}$ the solution to $\widetilde{K}(\alpha) = \max_{\bx\in \cX}\{\widetilde{g}(\bx,\alpha)\}$. We can compute the gradient, $\widetilde{K}'(\alpha)$, using the Danskin's theorem, as follows:
\begin{align}
    |\widetilde{K}'(\alpha)| &= \left(\sum_{s\in\cL}\sum_{\tau \in \Delta(s)} \exp\left(\frac{U(\tau;\widetilde{\bx})}{\mu}\right)\right) \geq  \sum_{s\in\cL}\sum_{\tau \in \Delta(s)} \exp\left(\frac{U(\tau;U^x \bbe)}{\mu}\right)\myeq{}\lambda, \nonumber 
\end{align}
Together with \eqref{eq:lemma59-eq1} we obtain the bound:
\begin{equation}
\label{eq:th8-eq2}
0\leq \delta^*-\overline{\delta}\leq \frac{\epsilon'(\cH+\cC)}{\lambda}.    
\end{equation}
Given the above bound on $\delta^*-\overline{\delta}$, we are now going to bound $\max_{\bx\in \cX}\widetilde{\cF}(\bx) - \widetilde{\cF}(\overline{\bx})$.
From the inequality $|\widetilde{g}(\overline{\bx},\overline{\delta})|\leq \epsilon' \cC$ claimed above (Equation~\ref{noname.1}), we have: 
\[
    \left| \sum_{s\in\cL}\sum_{\tau \in \Delta(s)}  r^l(s, \overline{x}_s) \exp\left(\frac{U(\tau;\overline{\bx})}{\mu}\right)  - \overline{\delta} \left(\sum_{s\in\cL}\sum_{\tau \in \Delta(s)} \exp\left(\frac{U(\tau;\overline{\bx})}{\mu}\right)\right)\right| \leq  \epsilon' \cC
\]
which implies
\[
\left|\widetilde{\cF}(\overline{\bx}) - \overline{\delta}\right| \leq \frac{\epsilon' \cC}{\sum_{s\in\cL}\sum_{\tau \in \Delta(s)} \exp\left(\frac{U(\tau;\overline{\bx})}{\mu}\right)} \leq \frac{\epsilon'
\cC}{\lambda}\implies \widetilde{\cF}(\overline{\bx})  \geq \overline{\delta} -\frac{\epsilon'
\cC}{\lambda}
\]
As a result, we obtain the following bounds:
\[
\delta^* = \max_{\bx\in \cX}\widetilde{\cF}(\bx)\geq \widetilde{\cF}(\overline{\bx})\geq \overline{\delta}  -\frac{\epsilon'
\cC}{\lambda} \implies \max_{\bx\in \cX}\widetilde{\cF}(\bx) - \widetilde{\cF}(\overline{\bx}) \leq \delta^* -\overline{\delta}+\frac{\epsilon'\cC}{\lambda} \stackrel{}{\leq} \frac{\epsilon'(\cH+2\cC)}{\lambda}.
\]
which concludes our proof.
\end{proof}
Finally, we present our theoretical bound with respect to our original problem \eqref{prob:main-homo}. Essentially, Theorem~\ref{th:overal-bound} is a direct result of Theorem \ref{thm.restricted.relation.2} and Theorem~\ref{prop.temp} with the  constant:
{\footnotesize\[
  \cU = \frac{\cH+2\cC}{(\min_{\bx\in\cX} \widetilde{\cF}(\bx)+\kappa) \sum_{s\in\cL}\sum_{\tau \in \Delta(s)} \exp\left(\frac{U(\tau;U^x\bbe)}{\mu}\right)}.
\]}
\begin{theorem}
\label{th:overal-bound}
If we run binary search to solve \eqref{prob:main-homo} with respect to the modified network $\cG'$ and obtain $(\overline{\bx},\overline{\delta})$, the following performance bound is guaranteed:
\[
\cF^l(\overline{\bx}) \geq \frac{1}{\overline{\eta}} \max_{\bx} \{\cF^l(\bx)\} -\kappa \frac{\overline{\eta}-1}{\overline{\eta}}, 
\]
where 
$
\overline{\eta} = (1+\beta_1)(1+\beta_2)\Big(1+ \epsilon' \cU
\Big), 
$ 
where $\cU$ is a constant independent of $\epsilon', \beta_1,\beta_2$.
\end{theorem}

Before  concluding, we give in Algorithm \ref{algo:algo1} the main steps of the approximation scheme. 
\begin{algorithm}[htb]
\SetAlgoLined
\DontPrintSemicolon
\caption{\textit{\textit{Solving \eqref{prob:main-homo} through the restricted interdiction problem \eqref{prob:approx-main}}}}\label{algo:main-approx-subprob} \label{algo:algo1}
\tcc{Solving \eqref{prob:approx-main} by approximately solving \eqref{eq:app-sub}}
- Create $\cG'$ by raising the travelling cost between any pair of nodes in $\cL$ \;
- Solve \eqref{prob:main-homo} using Algorithm \ref{alg:first} with the modified graph $\cG'$\;
\tcc{Improvement}
Recover the original graph $\cG$ and re-solve \eqref{prob:main-homo} to improve the
solution obtained from the above step.
\end{algorithm}

\begin{remark}
If $\beta_1$, $\beta_2$ and $\epsilon'$ are small, then $\overline{\eta}$ would be close to 1 and  $\overline{\bx}$ would be close to an optimal solution to the original interdiction problem. {Note that $\beta_1$, $\beta_2$ would be small in a real situation where the costs of traveling between critical nodes are expensive to the adversary (e.g., critical nodes are far away from each other). Moreover, if main paths for traveling between critical nodes can be well identified, it would be easy to raise the costs of these paths to make $\epsilon'$ arbitrarily small.} 
\end{remark}

\section{Numerical Experiments}

To illustrate the efficacy of our proposed algorithms, we perform experiments on synthetic data. 
\subsection{Experimental Settings}
\paragraph{Data generation:} We generate  random graphs (cycle-free)  with $\lvert \cS \rvert$ vertices and edge probability $p$. 
We randomly choose $\lvert \cL \rvert$ vertices (except source and destination) as the critical nodes that can be attacked. We set $\lvert \cL \rvert = 0.8 \times \lvert \cS \rvert$. In addition, the defender weights $\{(w^l_j, t^l_j) | j \in [\lvert \cL \rvert]\}$ are generated uniformly at random from the interval $[0, 1]$ and the adversary weights $\{(w^f_j, t^f_j) | j \in [\lvert \cS \rvert]\}$ are generated at random from the interval $[-1, 0]$. 
\paragraph{Baseline:} We approximate the sums over exponentially many paths in Equation~\ref{eq:trajectory_formulation} by sampling paths from the network and run gradient descent on this expression to estimate the optimal decision variable. Details on how paths are sampled are given in the appendix.
We use 1000 randomly sampled paths to estimate the baseline objective. To ensure fairness, all algorithms were run with the same number of epochs. Additionally, since Algorithm~\ref{algo:main-exact} has two gradient loops, we make sure that the sum of loops in Algorithm~\ref{algo:main-exact} is same as Baseline and Algorithm \ref{alg:first}. The \emph{run time} for all our methods for the largest graph size we consider (size $\lvert \cS \rvert = 100$) is \emph{under five minutes}. All our experiments were run on a 2.1 GHz CPU with 128GB RAM.
\subsection{Numerical Results} 
As shown in Table \ref{tab:comparisons}, we note that the Baseline struggles with solution quality compared to our proposed approaches. Moreover, the benefits of the derived guarantees of Algorithm~\ref{algo:main-exact} are also observed in our experiments.
\begin{table*}[htp]
    \centering
    \begin{tabular}{cccccc}
    \toprule
    & \multicolumn{5}{c}{Number of nodes ($\lvert \cS \rvert$)}\\
    Method & 20 & 40 & 60 & 80 & 100\\
    \midrule
       Baseline & 93.44 $\pm$ 3.12  & 94.38 $\pm$ 3.27 & 96.88 $\pm$ 2.18 & 94.72 $\pm$ 2.25 & 96.76 $\pm$ 2.53\\
        Algorithm \ref{alg:first} & 94.36 $\pm$ 2.69 & 98.10 $\pm$ 1.98 & 99.76 $\pm$ 0.16 & 99.48 $\pm$ 0.45 & 99.99 $\pm$ 0.00\\
    \bottomrule
    \end{tabular}
    \caption{Objective values of the optimal solutions obtained from various methods as a percentage of the optimal objective value obtained from Algorithm~\ref{algo:main-exact}. We used $p=0.8, \mu=2,
    \lvert \cL \rvert = 0.8 \lvert \cS \rvert$. 20 datasets were randomly generated for each setting and the mean and standard deviation are reported.}
    \label{tab:comparisons}
\end{table*}
\section{Conclusion}
Network interdiction game problems present a set of challenges that appear intractable to start with. In this work, we address some of these challenges by providing methods that solve a class of network interdiction problems with approximation guarantees. The quality of the approximation guarantee is problem dependent, but our methods empirically perform better than baselines over many randomly sampled problems. We are also the first to study the dynamic Quantal Response model in the type of network interdiction studied in \citep{fulkerson1977maximizing,israeli2002shortest}. We believe this modeling and methodology contribution provides suggestions for many possible future research directions, such as identifying properties of graphs for which our approximation guarantee is good, or studying a variant where the adversary's objective is to maximize a flow through the network.
%
%
%
%
%

\bibliographystyle{ACM-Reference-Format}
\bibliography{refs}
\newpage
\appendix

\section{Detailed Proofs}
\subsection{\textbf{Proof of Proposition~\ref{prop.1}.}}
\begin{proof}
We re-write the adversary's expected utility as follows:
\begin{align}
\cE^f &= \sum_{\tau \in \Omega} U(\tau)\frac{\exp\left(\frac{U(\tau)}{\mu}\right)}{\sum_{\tau'}\exp\left(\frac{U(\tau')}{\mu}\right)} \nonumber \\
&= \sum_{\tau \in \Omega^*}U(\tau)\frac{\exp\left(\frac{U(\tau)}{\mu}\right)}{\sum\limits_{\tau'\in \Omega^*}\exp\left(\frac{U(\tau')}{\mu}\right)+ \sum\limits_{\tau'\in \Omega\backslash\Omega^*}\exp\left(\frac{U(\tau')}{\mu}\right)} \nonumber \\
&+ \sum_{\tau \in \Omega\backslash\Omega^*}U(\tau)\frac{\exp\left(\frac{U(\tau)}{\mu}\right)}{\sum\limits_{\tau'\in \Omega^*}\exp\left(\frac{U(\tau')}{\mu}\right)+ \sum\limits_{\tau'\in \Omega\backslash\Omega^*}\exp\left(\frac{U(\tau')}{\mu}\right)}\nonumber\\
&= \frac{U(\tau^*)|\Omega^*|}{|\Omega^*|+ \sum\limits_{\tau'\in \Omega\backslash\Omega^*}\exp\left(\frac{U(\tau')-U(\tau^*)}{\mu}\right)}+ \sum\limits_{\tau \in \Omega\backslash\Omega^*}\frac{U(\tau)\exp\left(\frac{U(\tau)-U(\tau^*)}{\mu}\right)}{|\Omega^*|+ \sum\limits_{\tau'\in \Omega\backslash\Omega^*}\exp\left(\frac{U(\tau')-U(\tau^*)}{\mu}\right)}\label{eq:p1-eq1}
\end{align}
Let us denote  $\cT = \sum_{\tau'\in \Omega\backslash\Omega^*}\exp\left(\frac{U(\tau')-U(\tau^*)}{\mu}\right)$ for notational simplicity. We see that: 
\begin{equation}\label{eq:p1-eq2}
    0\leq \cT \leq |\Omega \backslash \Omega^*|\exp\left(-\frac{\alpha}{\mu}\right)
\end{equation}
From \eqref{eq:p1-eq1}, we have:
\begin{align}
    |\cE^f - U(\tau^*)|&\leq \frac{\cT}{|\Omega^*|+\cT} + \left| \sum_{\tau \in \Omega\backslash\Omega^*}\frac{U(\tau)\exp\left(\frac{U(\tau)-U(\tau^*)}{\mu}\right)}{|\Omega^*|+ \cT}\right|\nonumber \\
    &\leq  \frac{\cT}{|\Omega^*|+\cT} + L^*\left|\frac{\cT}{|\Omega^*| + \cT}\right|\nonumber \\
    &=(L^*+1)\frac{1}{1+\frac{|\Omega^*|}{\cT}} \stackrel{(i)}{\leq} \frac{L^*+1}{1+\frac{|\Omega^*|}{|\Omega\backslash\Omega^*|}\exp\left(\frac{\alpha}{\mu}\right)},\nonumber
\end{align}
where $(i)$ is due to \eqref{eq:p1-eq2}. We obtain the desired inequality. The limit $\lim_{\mu \rightarrow 0} \cE^f = U(\tau^*)$ is just a direct result of this inequality, concluding our proof. 
\end{proof}

\subsection{\textbf{Proof of Proposition~\ref{prop.2}}}

\begin{proof}
For any $n\in \bbN$, let us consider $\bM^n = \underbrace{\bM \times \bM \times \ldots \bM}_\text{$n$ times}$ with entries
\[
\bM^n_{ss'} = \sum_{\substack{(s_0,s_1,\ldots,s_n)\in \cS^{n+1} \\ s_0 = s, s_n = s'}} \left(\prod_{i=0}^{n-1}M_{s_is_{i+1}}  \right)
\]
Recall that $M_{ss'} = \exp\left(\frac{v(s)}{\mu}\right)$ if $s'\in N(s)$ and $M_{ss'} = 0$ otherwise. We see that if $n>|\cS|$, then  for any sequence $(s_0,s_1,\ldots,s_n)$ there is at least a pair $s_j = s_k$, $0\leq j,k\leq n $. Since  the network is cycle-free, there is at least  a pair $(s_j,s_{j+1})$ such that $M_{s_js_{j+1}} = 0$, leading to the fact that $\bM^n_{ss'} = \sum_{\substack{(s_0,s_1,\ldots,s_n)\in \cS^{n+1} \\ s_0 = s, s_n = s'}} \left(\prod_{i=0}^{n-1}M_{s_is_{i+1}}  \right) = 0$ for any $s,s'\in \cS$. Thus, if $n>|\cS|$ we have $\bM^n = \textbf{0}$. We select $n>|\cS|$ and  write
\[
 (\bI-\bM)\left(\sum_{t=0}^{n-1} \bM^{t} \right) = (\bI - \bM^n) =\bI,
\]
which implies $\det(\bI-\bM) \neq 0$, or equivalently, $\bI-\bM$ is invertible as desired.  
\end{proof}

\subsection{\textbf{Proof of Lemma~\ref{lemma.1}}}
\begin{proof}
Remind that we have $\beta_1$ and $\beta_2$ defined as follows:
\begin{equation*}
\begin{aligned}
    \beta_1 &= \max_{\bx}\max_{s\in \cL} \left\{\frac{\sum_{\tau \in \Delta^+(s)}  \exp\left(\frac{U(\tau;\bx)}{\mu}\right) }{\sum_{\tau \in \Delta(s)}  \exp\left(\frac{U(\tau;\bx)}{\mu}\right)}\right\}, &\beta_2 = \max_{\bx} \left\{\frac{\sum_{\tau \in \bigcup_s\{\Delta^+(s)\}}  \exp\left(\frac{U(\tau;\bx)}{\mu}\right) }{\sum_{\tau \in \bigcup_s\{\Delta(s)\}}  \exp\left(\frac{U(\tau;\bx)}{\mu}\right)}\right\},
\end{aligned}
\end{equation*}
For any defender strategy $x\in \cX$, we can re-write the defender's expected utility as follows:
\begin{align}
    \cF^l(\bx) &= \frac{\sum_{s\in\cL} r^l(s,x_s)\Bigg(\sum_{\tau\in \Delta(s)} \exp\left(\frac{U(\tau;\bx)}{\mu}\right) + \sum_{\tau\in \Delta^+(s)} \exp\left(\frac{U(\tau;\bx)}{\mu}\right)\Bigg)}{\sum_{\tau\in \bigcup_s\{\Delta(s)\}} \exp\left(\frac{U(\tau;\bx)}{\mu}\right) + \sum_{\tau\in \bigcup_s\{\Delta^+(s)\}} \exp\left(\frac{U(\tau;\bx)}{\mu}\right)} \myeq{}  \frac{\cU}{\cV} 
\end{align}
According to the definition of $\beta_1$, we obtain:
\begin{align*}
    & \sum_{\tau\in \Delta^+(s)} \exp\left(\frac{U(\tau;\bx)}{\mu}\right) \leq \beta_1\sum_{\tau\in \Delta(s)}  \exp\left(\frac{U(\tau;\bx)}{\mu}\right)\\
    \implies &\cU \leq \sum_{s\in\cL} r^l(s,x_s)\left(\sum_{\tau\in \Delta(s)} (1+\beta_1) \exp\left(\frac{U(\tau;\bx)}{\mu}\right)\right)
\end{align*}
In addition, we have:
    $\cV \geq \sum_{\tau\in \bigcup_s\{\Delta(s)\}} \exp\left(\frac{U(\tau;\bx)}{\mu}\right)$. 
As a result, we obtain the following inequality:
\begin{align*}
    & \cF^l(\bx) \le  \frac{\sum_{s\in\cL} r^l(s,x_s)\left(\sum_{\tau\in \Delta(s)} (1+\beta_1) \exp\left(\frac{U(\tau;\bx)}{\mu}\right)\right)}{\sum_{\tau\in \bigcup_s\{\Delta(s)\}} \exp\left(\frac{U(\tau;\bx)}{\mu}\right) }= ({1+\beta_1})\widetilde{\cF}(\bx)  &(*)
\end{align*}
On the other hand, according to the definition of $\beta_2$, we obtain: 
\begin{align*}
    & \sum_{\tau\in \bigcup_s\{\Delta^+(s)\}} \exp\left(\frac{U(\tau;\bx)}{\mu}\right)\leq \beta_2\sum_{\tau\in \bigcup_s\{\Delta(s)\}} \exp\left(\frac{U(\tau;\bx)}{\mu}\right) \\
    \implies &\cV \leq \sum_{\tau\in \bigcup_s\{\Delta(s)\}} \exp\left(\frac{U(\tau;\bx)}{\mu}\right) (1+\beta_2) 
\end{align*}
In addition, we have: $\cU \geq \sum_{s\in\cL} r^l(s,x_s)\left(\sum_{\tau\in \Delta(s)} \exp\left(\frac{U(\tau;\bx)}{\mu}\right)\right)$. As a result, we obtain:
\begin{align*}
    &\cF^l(\bx) \ge  \frac{\sum_{s\in\cL} r^l(s,x_s)\left(\sum_{\tau\in \Delta(s)} \exp\left(\frac{U(\tau;\bx)}{\mu}\right)\right)}{\sum_{\tau\in \bigcup_s\{\Delta(s)\}} \exp\left(\frac{U(\tau;\bx)}{\mu}\right) (1+\beta_2) }= \frac{1}{1+\beta_2}\widetilde{\cF}(\bx) & (**)
\end{align*}
The combination of (*) and (**) concludes our proof.
\end{proof}
\subsection{\textbf{Proof of Lemma~\ref{lm:choose-kappa}}}
\begin{proof}
We reuse the definitions of $\cU$ and $\cV$ as in the proof of Lemma~\ref{lemma.1}. Similar to proof of Lemma~\ref{lemma.1}, according to the definition of $\beta_1$, we obtain:
\begin{align*}
    & \sum_{\tau\in \Delta^+(s)} \exp\left(\frac{U(\tau;\bx)}{\mu}\right) \leq \beta_1\sum_{\tau\in \Delta(s)}  \exp\left(\frac{U(\tau;\bx)}{\mu}\right)\\
\end{align*}
Besides, according to the definition $\kappa =\sum_{s\in \cL}\max_{\bx}\left|r^l(s,x_s) \right|$, we have $r^l(s,x_s)+\kappa\geq 0 $ for any $\bx\in\cX$. Thus, we can write: 
\begin{align*}
     \cU + \kappa \cV &= \sum_{s\in\cL} \left(r^l(s,x_s)+\kappa\right)\left(\sum_{\tau\in \Delta(s) \bigcup \Delta^+(s)} \exp\left(\frac{U(\tau;\bx)}{\mu}\right)\right) \\
    &\leq  \sum_{s\in\cL} \left(r^l(s,x_s) + \kappa\right)\left(\sum_{\tau\in \Delta(s)} (1+\beta_1) \exp\left(\frac{U(\tau;\bx)}{\mu}\right)\right)
\end{align*}
In addition, we have:
    $\cV \geq \sum_{\tau\in \bigcup_s\{\Delta(s)\}} \exp\left(\frac{U(\tau;\bx)}{\mu}\right)$. 
As a result, we obtain the following inequality:
\begin{align}\nonumber
     \cF^l(\bx) + \kappa = \frac{\cU + \kappa\cV}{\cV}&\le  \frac{\sum_{s\in\cL} \left(r^l(s,x_s) + \kappa\right)\left(\sum_{\tau\in \Delta(s)} (1+\beta_1) \exp\left(\frac{U(\tau;\bx)}{\mu}\right)\right)}{\sum_{\tau\in \bigcup_s\{\Delta(s)\}} \exp\left(\frac{U(\tau;\bx)}{\mu}\right) }\\
     &= ({1+\beta_1})\left(\widetilde{\cF}(\bx) + \kappa\right) \label{eq:b1}
\end{align}
On the other hand, from the way we select $\kappa$, we have:
\begin{align}
    &\kappa \geq \sum_{s\in\cL} |r^l(s,x_s)| \geq \sum_{s\in\cL} |r^l(s,x_s)| \frac{\sum_{\tau\in \Delta(s)} \exp\left(\frac{U(\tau;\bx)}{\mu}\right)}{\sum_{\tau\in \bigcup_{s'}\{\Delta(s')\}} \exp\left(\frac{U(\tau;\bx)}{\mu}\right)} \nonumber\\
    \implies &\sum_{s\in\cL} \left(r^l(s,x_s)+\kappa\right) \left(\sum_{\tau\in \Delta(s)} \exp\left(\frac{U(\tau;\bx))}{\mu}\right) \right)\geq 0\nonumber
\end{align}
Besides, according to the definition of $\beta_2$, we obtain: 
\begin{align*}
    & \sum_{\tau\in \bigcup_s\{\Delta^+(s)\}} \exp\left(\frac{U(\tau;\bx)}{\mu}\right)\leq \beta_2\sum_{\tau\in \bigcup_s\{\Delta(s)\}} \exp\left(\frac{U(\tau;\bx)}{\mu}\right) \\
    \implies &\cV \leq \sum_{\tau\in \bigcup_s\{\Delta(s)\}} \exp\left(\frac{U(\tau;\bx)}{\mu}\right) (1+\beta_2) 
\end{align*}
As a reulst, we obtain the following inequalities,
\begin{align}
\cF^l(\bx) +\kappa &\geq { \frac{\sum_{s\in\cL} (r^l(s,x_s)+\kappa) \left(\sum_{\tau\in \Delta(s)} \exp\left(\frac{U(\tau;\bx))}{\mu}\right) \right)}
    {\cV}}\nonumber \\
    &\geq { \frac{\sum_{s\in\cL} (r^l(s,x_s)+\kappa) \left(\sum_{\tau\in \Delta(s)} \exp\left(\frac{U(\tau;\bx))}{\mu}\right) \right)}
    {\sum_{\tau\in \bigcup_s\{\Delta(s)\}} \exp\left(\frac{U(\tau;\bx)}{\mu}\right) (1+\beta_2)}} = \frac{1}{1+\beta_2} (\widetilde{\cF}(\bx)+\kappa)\label{eq:b2}.
\end{align}
Combining \eqref{eq:b1} and \eqref{eq:b2} gives us the desired inequalities. 
\end{proof}

\subsection{Proof of Lemma~\ref{lemma.temp}}
\begin{proof}
Observing that $\{\tau; \tau\in \Omega, \tau \ni s \} = \Delta^+(s) \cup \Delta(s)$ and $\Delta^+(s), \Delta(s)$ are disjoint, we can decompose $g(\bx,\delta|\cG')$ into two separate terms, as follows:
\begin{align}
g(\bx,\delta|\cG') &= \sum_{s\in\cL}\sum_{\substack{\tau\in \Omega \\ \tau \ni s}}  r^f_s(x_s) \exp\left(\frac{U(\tau;\bx)}{\mu}\right) - \delta \left(\sum_{\tau'\in\Omega} \exp\left(\frac{U(\tau';\bx)}{\mu}\right)\right) \nonumber \\
&= \widetilde{g}(\bx, \delta) + \cT(\bx, \delta),
\end{align}
where the second term: 
\begin{align*}
    \cT(\bx,\delta) = \sum_{s\in\cL}r^f_s(x_s) \sum_{\substack{\tau\in \Delta^+(s)}}   \exp\left(\frac{U(\tau;\bx)}{\mu}\right)
    - \delta \left(\sum_{\tau\in\bigcup_{s\in\cL} \Delta^+(s) } \exp\left(\frac{U(\tau;\bx)}{\mu}\right)\right).
\end{align*}
Moreover, remind that we have the definition of $\rho^s$ and $\rho$:
\begin{align}\nonumber
    &\rho^s = \sum_{\substack{\tau\in \Omega,\tau\ni s}}   \exp\left(\frac{U(\tau;L^x \bbe)}{\mu}\right),\;\forall s\in \cL 
    & \rho = \sum_{\substack{\tau\in \Omega}}   \exp\left(\frac{U(\tau;L^x \bbe)}{\mu}\right)
\end{align}
According to conditions in Equation~\eqref{eq:eps-1} and \eqref{eq:eps-2}, we obtain:
\begin{align}
     &\max_{\bx}\max_{s\in \cL} \left\{\frac{\sum_{\tau \in \Delta^+(s)}  \exp\left(\frac{U(\tau;\bx)}{\mu}\right) }{\sum_{\tau \in \Omega, \tau \ni s}  \exp\left(\frac{U(\tau;\bx)}{\mu}\right)}\right\} \leq \epsilon' \implies \sum_{\substack{\tau\in \Delta^+(s)}}   \exp\left(\frac{U(\tau;\bx)}{\mu}\right)\le \epsilon' \rho^s,\;\forall \bx\in \cX, s\in \cL \label{eq:prop58-eq3} \\
    &\max_{\bx}\left\{\frac{\sum_{\tau \in \bigcup_s\{\Delta^+(s)\}}  \exp\left(\frac{U(\tau;\bx)}{\mu}\right) }{\sum_{\tau \in\Omega}  \exp\left(\frac{U(\tau;\bx)}{\mu}\right)}\right\} \leq \epsilon'\implies \sum_{\substack{\tau \in \bigcup_s\{\Delta^+(s)\}}}   \exp\left(\frac{U(\tau;\bx)}{\mu}\right) \le \epsilon' \rho,\;\forall \bx\in \cX. \label{eq:prop58-eq4}  
\end{align}
By using these inequalities, we have: 
\begin{align}
    |\cT(\bx,\delta)| &\leq \sum_{s\in\cL}  |r^f_s(x_s)|\epsilon'\rho^s + \delta \epsilon' \rho \leq \epsilon'\left(\kappa \max_s\{\rho^s\} + \delta \rho\right). 
\end{align}
which concludes our proof.
\end{proof}
\subsection{Proof of Lemma~\ref{prop:gap-g-tilde}}
\begin{proof}

We reuse the decomposition of $g(\bx,\delta|\cG')$ as described in Lemma~\ref{lemma.temp}. By taking the derivative of  $\cT(\bx,\delta)$ w.r.t $x_j$, $j\in\cL$, we obtain: 
\begin{align*}
     \frac{\partial \cT(\bx,\delta)}{\partial x_j} &= w^f_j \sum_{\substack{\tau\in \Delta^+(j)}}   \exp\left(\frac{U(\tau;\bx)}{\mu}\right) +  \frac{1}{\mu}\sum_{s\in\cL}r^f_s(x_s) \sum_{\substack{\tau\in \Delta^+(s),\tau \ni j}} w^f_j  \exp\left(\frac{U(\tau;\bx)}{\mu}\right)\nonumber \\
     &- \frac{\delta}{\mu} \left(\sum_{\tau\in\bigcup_{s\in\cL} \Delta^+(s),\tau\ni j } w^f_j \exp\left(\frac{U(\tau;\bx)}{\mu}\right)\right).
\end{align*}
Thus, using \eqref{eq:prop58-eq3} and \eqref{eq:prop58-eq4}, we get the inequality:
\begin{equation}
\label{eq:b3}
\left|\frac{\partial \cT(\bx,\delta)}{\partial x_j}\right| \leq  \epsilon' \left( |w^f_j|\rho^j + \frac{\kappa}{\mu} \max_{s}\{\rho^s\} +\frac{\delta}{\mu}  |w^f_j| \rho\right).
\end{equation}
Now, let us consider the problem $\max_{\bx}\{g(\bx,\delta|\cG')|\; \bx\in \cX\}$. We form its Lagrange dual as: 
\begin{align*}
    L^{\cG'}(\bx, \bgamma,\bbeta)   = {g}(\bx, \delta|\cG') - \sum_{k\in[K]}\gamma_k \left(\sum_{s\in\cL_k}x_s
    -  M_k\right) - \sum_{s}\eta^1_s (x_s - U^x) + \sum_{s}\eta^2_s (x_s-L^x).
\end{align*}
If $\overline{\bx}$ is a stationary point of $\max_{\bx}\{g(\bx,\delta|\cG')|\; \bx\in \cX\}$, then the KKT conditions imply that there are $\bgamma^*,\bbeta^{1*}$, $\bbeta^{2*}$ such that:
\[
 \frac{{\partial g}(\overline{\bx}, \delta|\cG')}{\partial x_j} - \sum_{k\in[K]}\gamma^*_k\bbI[j\in \cL_k]  - \eta^{1*}_j  + \eta^{2*}_j = 0, 
\]
where $\bbI[\cdot]$ is the indicator function. Note that $\frac{\partial {g}(\overline{\bx}, \delta|\cG')}{\partial x_j}  = \frac{\partial \widetilde{g}(\overline{\bx}, \delta)}{\partial x_j} + \frac{\partial\cT(\bx,\delta) }{\partial x_j} $. Thus, from \eqref{eq:b3} we have:
\begin{align}\nonumber
& \left|\frac{\partial \widetilde{g}(\overline{\bx}, \delta)}{\partial x_j} - \sum_{k\in[K]}\gamma^*_k \bbI[j\in \cL_k]  - \eta^{1*}_j  + \eta^{2*}_j\right| \\
 &\leq 
 \left|\frac{{\partial g}(\overline{\bx}, \delta|\cG')}{\partial x_j} - \sum_{k\in [K]}\gamma^*_k \bbI[j\in \cL_k]  - \eta^{1*}_j  + \eta^{2*}_j\right| + \left|\frac{\partial\cT(\bx,\delta) }{\partial x_j}\right| 
 \nonumber \\&\leq  \epsilon' \left( |w^f_j|\rho^j + \frac{\kappa}{\mu} \max_{s}\{\rho^s\} +\frac{\delta}{\mu}  |w^f_j| \rho\right).\label{eq:bound-grad-gt}  
\end{align}
Let us now  define a function $\widehat{G}(\bx,\delta)$ as follows: $$
\widehat{G}(\bx,\delta) = \widetilde{g}(\bx,\delta) +\sum_{s\in \cL} \alpha_s \exp\left(\frac{w^f_sx_s+t^f_s}{\mu}\right),$$ where $\alpha_s$, $\forall s\in \cL$, are chosen as follows:
\[
\alpha_s = -\frac{\frac{\partial \widetilde{g}(\overline{\bx},\delta)}{\partial x_s} - \sum_{k\in[K]}\gamma^*_k \bbI[j\in \cL_k]  - \eta^{1*}_j  + \eta^{2*}_j}{\frac{w^f_s}{\mu}\exp\left(\frac{w^f_sx_s+t^f_s}{\mu}\right)}.
\]
Given $\alpha_s$ defined as above, we obtain the following equations:
\begin{equation}
\label{eq:lm4-eq11}
\frac{\partial \widehat{G}(\overline{\bx},\delta)}{\partial x_j}  - \sum_{k\in[K]}\gamma^*_k \bbI[j\in \cL_k]  - \eta^{1*}_j  + \eta^{2*}_j= 0,~\forall j\in \cL.
\end{equation}
In the following, we first attempt to bound the gap $\left|\widehat{G}(\bx, \delta) - \widetilde{g}(\bx,\delta)\right|$ for every defender strategy $\bx$. We then leverage this bound together with the unimodality of $\widehat{G}(\bx, \delta)$ to bound the gap $|\max_{\bx\in\cX}\widetilde{g}(\bx,\delta) -\widetilde{g}(\overline{\bx},\delta)|$. As we show later, the unimodality of $\widehat{G}(\bx, \delta)$ is proved based on Equation~\ref{eq:lm4-eq11} and the variable conversion trick used  in the proof of Theorem~\ref{prop:c1}. 

\paragraph{\textbf{Bounding the gap $\left|\widehat{G}(\bx, \delta) - \widetilde{g}(\bx,\delta)\right|$.}}By taking the derivative of $\widehat{G}(\bx,\delta)$ w.r.t. $x_j$, we get:
\[
\frac{\partial  \widehat{G}(\bx, \delta)}{\partial x_j} = \frac{ \widetilde{g}(\bx,\delta)}{\partial x_j} + \frac{\alpha_j w^f_j}{\mu} \exp\left(\frac{w^f_j+t^f_j}{\mu}\right). 
\]
Combining with \eqref{eq:bound-grad-gt}-\eqref{eq:lm4-eq11} we can bound $\alpha_j$, $\forall j\in \cL$,  as follows:
\begin{align}
    &\left|\frac{\alpha_j w^f_j}{\mu} \exp\left(\frac{w^f_j\overline{x}_j+t^f_j}{\mu}\right)\right| = \left|\frac{\partial  \widehat{G}(\overline{\bx}, \delta)}{\partial x_j} - \frac{ \widetilde{g}(\overline{\bx},\delta)}{\partial x_j}\right|\nonumber \\ 
    &\leq \left| \frac{\partial \widehat{G}(\overline{\bx},\delta)}{\partial x_j}  - \sum_{k\in[K]}\gamma^*_k \bbI[j\in \cL_k]  - \eta^{1*}_j  + \eta^{2*}_j\right| + \left| \frac{\partial \widetilde{g}(\overline{\bx},\delta)}{\partial x_j}  - \sum_{k\in[K]}\gamma^*_k \bbI[j\in \cL_k]  - \eta^{1*}_j  + \eta^{2*}_j\right|\nonumber\\
    &\leq \epsilon' \left( |w^f_j|\rho^j + \frac{\kappa}{\mu} \max_{s}\{\rho^s\} +\frac{\delta}{\mu}  |w^f_j| \rho\right)
\end{align}
which implies: 
\begin{equation}
\label{eq:lm4-eq12}
    \left|\alpha_j\exp\left(\frac{w^f_j\overline{x}_j+t^f_j}{\mu}\right)\right| \leq  {\epsilon'\mu} \left( \rho^j + \frac{\kappa}{\mu|w^f_j|} \max_{s}\{\rho^s\} +\frac{\delta}{\mu}   \rho\right).
\end{equation}
Combining the above inequality with the definition of $\widehat{G}(\bx, \delta)$ we obtain the following inequality for all defender strategy $\bx\in\cX$:
\begin{align}
    \left|\widehat{G}(\bx, \delta) - \widetilde{g}(\bx,\delta)\right| &\leq\sum_{j\in \cL} \left|\alpha_j\exp\left(\frac{(w^f_j{x}_j+t^f_j)}{\mu}\right)\right|\nonumber\\
    &\leq \sum_{j\in \cL} \left|\alpha_j\exp\left(\frac{w^f_j\overline{x}_j+t^f_j}{\mu}\right) \exp\left(\frac{w^f_j(x_j-\overline{x}_j)}{\mu}\right)\right|\nonumber\\
    &\leq   {\epsilon'} \left( \mu\rho^j + \frac{\kappa}{|w^f_j|} \max_{s}\{\rho^s\} +{\delta}   \rho\right)\exp\left(\frac{w^f_j(L^x-U^x)}{\mu}\right).\nonumber
\end{align}
\paragraph{\textbf{Bounding $|\max\limits_{\bx\in\cX}\widetilde{g}(\bx,\delta) -\widetilde{g}(\overline{\bx},\delta)|$.}}Let us define $\cH =   2\left( \mu\rho^j + \frac{\kappa}{|w^f_j|} \max_{s}\{\rho^s\} +{\delta}   \rho\right)\exp\left(\frac{w^f_j(L^x-U^x)}{\mu}\right) $ for notational simplicity. We have the following remarks:
\begin{itemize}
    \item[(i)] From \eqref{eq:lm4-eq11}, we see that $\overline{\bx}$ is a stationary point of the maximization problem $\max_{\bx\in \cX} \{\widehat{G}(\bx,\delta)\}$  and $\bgamma^*,\bbeta^{1*}$, $\bbeta^{2*}$ are the corresponding KKT multipliers.
    \item[(ii)] With the change of variables used  in the proof of Theorem~\ref{prop:c1}, the function $\widehat{G}(\bx,\delta)$ becomes $\widetilde{g}(\widehat{\bx}(\by),\delta) + \sum_{s\in \cL} \alpha_s y_s$, thus  $\widehat{G}(\widehat{\bx}(\by),\delta)$ is strictly concave in $\by$. Similar to Theorem~\ref{prop:c1}, $\widehat{G}(\bx,\delta)$ is unimodal (i.e., any local optimum is a  global one). 
\end{itemize}
Thus, $\overline{\bx}$ is also an optimal solution to $\max_{\bx\in \cX} \widehat{G}(\bx,\delta)$. As a result, we now can bound the gap $|\max_{\bx\in\cX}\widetilde{g}(\bx,\delta) -\widetilde{g}(\overline{\bx},\delta)|$ as follows:
\begin{align}
    |\max_{\bx\in\cX}\widetilde{g}(\bx,\delta) -\widetilde{g}(\overline{\bx},\delta)| &\leq |\max_{\bx\in\cX}\widetilde{g}(\bx,\delta) -\widehat{G}(\overline{\bx},\delta)| + |\widehat{G}(\overline{\bx},\delta) - \widetilde{g}(\overline{\bx},\delta)| \nonumber \\
    &\leq |\max_{\bx\in\cX}\widetilde{g}(\bx,\delta) -\max_{\bx\in \cX}\widehat{G}({\bx},\delta)| +\frac{\epsilon' \cH}{2}.\label{eq:lm4-eq21}
\end{align}
We consider the following two cases:
\begin{itemize}
    \item If $\max_{\bx\in\cX}\widetilde{g}(\bx,\delta) \geq\max_{\bx\in \cX}\widehat{G}({\bx},\delta)$. Let  $\bx^*$ be optimal for $\max_{\bx\in\cX}\widetilde{g}(\bx,\delta)$, we have
    \begin{align}
        |\max_{\bx\in\cX}\widetilde{g}(\bx,\delta) -\max_{\bx\in \cX}\widehat{G}({\bx},\delta)|&= \widetilde{g}(\bx^*,\delta) - \max_{\bx\in \cX}\widehat{G}({\bx},\delta) \nonumber \\
        &\leq \widetilde{g}(\bx^*,\delta) - \widehat{G}({\bx^*},\delta)\nonumber \\
        &\leq \frac{\epsilon'\cH}{2} \label{eq:lm4-eq22}
    \end{align}
    \item If $\max_{\bx\in\cX}\widetilde{g}(\bx,\delta) <\max_{\bx\in \cX}\widehat{G}({\bx},\delta)$, then  we have
    \begin{align}
        |\max_{\bx\in\cX}\widetilde{g}(\bx,\delta) -\max_{\bx\in \cX}\widehat{G}({\bx},\delta)|&= \max_{\bx\in \cX}\widehat{G}({\bx},\delta)  - \max_{\bx\in \cX} \widetilde{g}(\bx,\delta) \nonumber \\
        &\leq \widehat{G}(\overline{\bx},\delta) - \widetilde{g}(\overline{\bx},\delta) \nonumber \\
        &\leq \frac{\epsilon'\cH}{2} \label{eq:lm4-eq23}
    \end{align}
\end{itemize}
Combine \eqref{eq:lm4-eq21}-\eqref{eq:lm4-eq22}-\eqref{eq:lm4-eq23} we obtain:
\[
|\max_{\bx\in\cX}\widetilde{g}(\bx,\delta) -\widetilde{g}(\overline{\bx},\delta)| \leq \epsilon'\cH
\]
which is the desired inequality, completing our proof.

\end{proof}
\section{Experimental Details}
To sample paths for the baseline, a resource allocation $\bx$ is assigned to $\cL$ and the follower is initially placed at the origin $s_o$. Its next node $s_1$ is sampled from the distribution $\pi^f(s|s_0, \bx) = \frac{\exp\left(\frac{v(s;\bx)}{\mu}\right)Z_s}{ \sum_{s'\in N(s_0)} \exp\left(\frac{v(s';\bx)}{\mu}\right)Z_{s'}}$ where $s \in \{$nodes having an edge to $s_0$\} and $N(s_0)$ is the set of outgoing nodes from $s_0$. Similarly $s_2, \ldots$ etc. are sampled till the destination $s_d$ is reached. This sampling is repeated 1000 times per iteration and then average is taken to get the objective. Based on the gradients, the resource allocation $\bx$ is updated which changes the transition probabilities and the process is repeated again until convergence. Ten different initializations of $\bx$ were taken and the seed with the lowest loss was reported.

\section{Zero-sum Game Model} 

\subsection{Problem Formulation}
In this section  we discuss  a zero-sum game model that is often used in adversarial settings, in which the aim of the defender is to minimize the expected utility of the adversary. {The adversary's expected utility can be computed as follows:
\[
\cE^f(\bx)  = \sum_{\tau \in \Omega} U(\tau;\bx)\frac{\exp\left(\frac{U(\tau;\bx)}{\mu}\right)}{\sum_{\tau'\in \Omega} \exp\left(\frac{U(\tau';\bx)}{\mu}\right)}
\]
The zero-sum game model can then be formulated as follows:
\begin{align}
     \min_{\bx\in\cX } &\quad  \cE^f(\bx)   \label{prob:zero-sum} \tag{\sf OPT-zerosum}
\end{align}
which is generally non-convex in $\bx$. Since it shares the same structure with the non-zero-sum game model considered in the main body of the paper, our approximation method based on the restricted problem still applies. Here, instead of directly solve the non-convex problem  \eqref{prob:zero-sum}, we propose to optimize the following log-sum objective, which is more tractable to handle
\[\Gamma(\bx) = \mu \log\Big(\sum\nolimits_{\tau\in\Omega} \exp\left(\frac{U(\tau;\bx)}{\mu}\right)\Big)
\]
} 
It can be seen that $\Gamma(\bx)$ has a log-sum-exp convex form of a geometric program, thus it is convex~\citep{boyd2004convex}. 
From the results in Section~\ref{sec:handling-exp-paths}, we further see that
$\Gamma(\bx)$ can be computed by solving  a system of linear equations, which can be done in poly-time. Thus, the optimization problem $\max_{\bx}\Gamma(\bx)$ can be solved in poly-time.  We discuss in the following a connection between  \eqref{prob:zero-sum}, the alternative formulation $\max_{\bx}\Gamma(\bx)$ and the a classical shortest-path network interdiction problem \citep{smith2020survey,israeli2002shortest}.
 To facilitate explanation of this point, let us consider the following shortest-path network interdiction problem:
\begin{align}
     \min_{\bx \in \cX }  &\quad \Big\{ T(\bx) = \max_{\tau \in \Omega} \quad  U(\tau;\bx) \Big\}. \label{prob:zero-sum-shortest-path} \tag{\sf OPT-shortest-path}
\end{align}
It is known that the above shortest-path network interdiction problem  can be formulated as a mixed-integer linear program and is NP-hard~\citep{israeli2002shortest}. We first  bound the gap between $\Gamma{(\bx)}$ and $T(\bx)$ for any $\bx\in \cX$ in Lemma \ref{lm:lm-zero-sum-1} below
\begin{lemma}
\label{lm:lm-zero-sum-1} For $\bx\in \cX$, let $\tau^* = \text{argmax}_{\tau \in \Omega}U(\tau; \bx)$ (i.e., the best trajectory which gives the highest adversary utility), $\Omega^* = \{\tau|\; U(\tau;\bx) = U(\tau^*;\bx)\}$ (i.e., the set of all trajectories with the same highest utility), and $\alpha = \{U(\tau^*;\bx) - \max_{\tau \in \Omega\backslash\Omega^*} U(\tau;\bx)\}$, then we have:
\[
 \left| \Gamma(\bx)- T(\bx)\right|\leq \mu \log \left(|\Omega^*|+\frac{|\Omega \backslash\Omega^*|}{\exp(\alpha/\mu)} \right).
\]
As a result, $\lim_{\mu\rightarrow 0 }\Gamma(\bx) =  U(\tau^*)$.
\end{lemma}
\begin{proof}
We can write: 
\begin{align}
\Gamma(\bx) &= \mu \log\left(\sum_{\tau\in \Omega} \exp\left(\frac{U(\tau;\bx)}{\mu}\right) \right)   \nonumber \\
&=  \mu \log\left(|\Omega^*|\exp\left(\frac{U(\tau^*;\bx)}{\mu}\right) +\sum_{\tau\in \Omega\backslash \Omega^*} \exp\left(\frac{U(\tau;\bx)}{\mu}\right) \right)   \nonumber \\
&\leq \mu \log\Bigg(|\Omega^*|\exp\left(\frac{U(\tau^*;\bx)}{\mu}\right)  +(|\Omega\backslash \Omega^*|) \exp\left(\frac{U(\tau^*;\bx)-\alpha}{\mu}\right) \Bigg)   \nonumber \\
&=U(\tau^*;\bx) + \mu \log\left(|\Omega^*|+ (|\Omega\backslash\Omega^*|)\exp\left(\frac{-\alpha}{\mu}\right)\right)\label{eq:eq1}
\end{align}
Moreover, we have $\Gamma(\bx) \geq U(\tau^*;\bx)$. Combine this with \eqref{eq:eq1} we obtain the desired inequality. The limit $\lim_{\mu \rightarrow 0} \Gamma(\bx) = U(\tau^*)$ is just a direct result of this equality, concluding our proof. 
\end{proof}
Combine Lemma \ref{lm:lm-zero-sum-1} with Proposition \ref{prop.1}, we obtain a bound for $|\cE^f(\bx)-\Gamma(\bx)|$ 
\begin{lemma}
\label{lm:lm-zero-sum-2}
Let $L^* = \max_{\tau\in \Omega,\bx\in \cX} |U(\tau;\bx)|$, we have 
\[
|\cE^f(\bx) - \Gamma(\bx)| \leq 
\frac{L^*+1}{1+\frac{|\Omega^*|}{|\Omega\backslash\Omega^*|}\exp\left(\frac{\alpha}{\mu}\right)} + \mu \log \left(|\Omega^*|+\frac{|\Omega \backslash\Omega^*|}{\exp\left(\frac{\alpha}{\mu}\right)} \right)
\]
\end{lemma}
We are now ready to assess the quality of a solution given by the alternative formula $\max_{\bx} \Gamma(\bx)$ and the zero-sum game ones \eqref{prob:zero-sum} and \eqref{prob:zero-sum-shortest-path}. Let $\Gamma^* = \max_{\bx}\Gamma(\bx)$, $\cE^*$, $T^*$ be the optimal value of \eqref{prob:zero-sum}, \eqref{prob:zero-sum-shortest-path}, and $\overline{\bx}$ be the optimal solution to $\max_{\bx}\Gamma(\bx)$. Given any $\bx\in\cX$, let: $$\alpha(\bx) = \max_{\tau\in\Omega} U(\tau;\bx) - \max \{U(\tau;\bx)|\; \tau\in\Omega,\;U(\tau;\bx) < \max_{\tau\in\Omega}U(\tau;\bx)\}$$ Intuitively, $\alpha(\bx)$ is the adversary loss in utility if the adversary chooses the second best trajectory instead of the optimal one. In addition, let $C(\bx)$ be the number of best paths in $\Omega$, that is, $C(\bx) = |\arg\!\max_{\tau \in \Omega} U(\tau;\bx)|$. We have the following  results bounding the gaps between the convex problem $\max_{\bx} \Gamma(\bx)$  and  two baselines, i.e., the classical shortest-path network interdiction and its bounded rational version,  as functions of $\mu$. The results imply that the optimal values and optimal solutions to $\max_{\bx} \Gamma(\bx)$ converge  to those of \eqref{prob:zero-sum-shortest-path} and \eqref{prob:zero-sum} when $\mu$ goes to zero. 
\begin{proposition}
\label{prop:zero-sum}
Let $\bx^* = \text{argmax}_{\bx} \Gamma(\bx)$ and  
\begin{align}
    \kappa_1(\mu) \!&=\!\max_{\bx}   \left\{\mu \log \left(|C({\bx})+\frac{|\Omega|-|C({\bx})}{\exp\left(\frac{\alpha({\bx})}{\mu}\right)} \right)\right\}\nonumber\\
    \kappa_2(\mu) \!&=\!  
\kappa_1(\mu) + \max_{\bx} \left\{\mu \frac{L^*+1}{1+\frac{|C({\bx})}{|\Omega|-|C({\bx})}\exp\left(\frac{\alpha({\bx})}{\mu}\right)}\right\}
\end{align}
The following results hold
\begin{itemize}
    \item[(i)]  $|\Gamma^* - T^*| \leq \kappa_1(\mu), \text{ and }
  |T({\bx^*}) - \max_{{\bx}} \{T(\bx)\}| \leq 2\kappa_1(\mu)$
  \item[(ii)]$|\Gamma^* - \cE^*| \leq \kappa_2(\mu), \text{ and }
  |\cE^f({\bx^*}) - \max_{{\bx}} \{\cE^f(\bx)\}| \leq 2\kappa_2(\mu)$
  \item[(iii)]$\lim\limits_{\mu\rightarrow 0} \kappa_1(\mu) = \kappa_2(\mu) = 0$
\end{itemize}
\end{proposition}
\begin{proof} For (i), we first note that $\Gamma(\bx) \leq T(\bx)$ for any $\bx\in \cX$. Thus $\Gamma^* \leq  T^*$. Let $\overline{\bx}$ be an optimal solution to \eqref{prob:zero-sum-shortest-path}, we  write:
\begin{align}
    |\Gamma^* - T^*| =  T^* - \Gamma^*&= T(\overline{\bx}) - \max_{\bx}  \Gamma(\bx)\nonumber \\
    &\leq T(\overline{\bx}) -  \Gamma(\overline{\bx})\nonumber \\
    &\stackrel{(a)}{\leq} \kappa_1(\mu).\nonumber
\end{align}
where $(a)$ is due to Proposition \ref{prop.1}. 
Moreover, considering the gap $|T({\bx^*}) - \max_{{\bx}} T(\bx)|$, we have the chain of inequalities
\begin{align}
    |T({\bx^*}) - \max_{{\bx}} \{T(\bx)\}|&\leq |T({\bx^*}) - \max_{{\bx}} \{\Gamma(\bx)\}|+ |\Gamma^* - T^*|\nonumber \\
    &= |T({\bx^*}) - \Gamma({\bx^*})|+ |\Gamma^* - T^*|\nonumber \\
    &\leq 2\kappa_1(\mu).\nonumber 
\end{align}
For (ii), let $\widehat{\bx}$ be an optimal solution to  \eqref{prob:zero-sum}. Similarly,  we can write 
\begin{align}
    |\Gamma^* - \cE^*| =  \cE^* - \Gamma^*&= \cE^f(\widehat{\bx}) - \max_{\bx}  \Gamma(\bx)\nonumber \\
    &\leq \cE^f(\widehat{\bx}) -  \Gamma(\widehat{\bx})\nonumber \\
    &\stackrel{(b)}{\leq} \kappa_2(\mu).
\end{align}
where $(b)$ is due to Lemma \ref{lm:lm-zero-sum-2}. Moreover, considering the gap $|\cE^f(\bx^*) - \max_{\bx}\cE^f(\bx)|$, we write
\begin{align}
    |\cE^f({\bx^*}) - \max_{{\bx}} \{\cE^f(\bx)\}|&\leq |\cE^f({\bx^*}) - \max_{{\bx}} \{\Gamma(\bx)\}|+ |\Gamma^* - \cE^*|\nonumber \\
    &= |\cE^f({\bx^*}) - \Gamma({\bx^*})|+ |\Gamma^* - \cE^*|\nonumber \\
    &\leq 2\kappa_2(\mu).\nonumber 
\end{align}
The limits  $\lim\limits_{\mu\rightarrow 0} \kappa_1(\mu) = \kappa_2(\mu) = 0$ are obviously verified, which concludes the proof.
\end{proof}
In fact, we can control the adversary's rationality by adjusting $\mu$, i.e., the adversary would be more rational as $\mu\!\rightarrow\! 0$ (perfectly rational if $\mu=0$), and be irrational as $\mu\!\rightarrow \!\infty$.
\subsection{Experiment Results for Zero-Sum Games}
\begin{table*}[t!]
    \centering
    \begin{tabular}{cccccc}
    \toprule
    & \multicolumn{5}{c}{Number of nodes ($\lvert \cS \rvert$)}\\
    Method & 20 & 40 & 60 & 80 & 100\\
    \midrule
       Baseline & 99.95 $\pm$ 0.00  & 99.88 $\pm$ 0.07 & 99.74 $\pm$ 0.18 & 99.69 $\pm$ 0.23 & 99.07 $\pm$ 0.52\\
    \bottomrule
    \end{tabular}
    \caption{Objective values of the optimal solutions obtained from the Baseline  as a percentage of the optimal objective value obtained using our approach for handling an exponential number of paths. We use $p=0.8, \mu=2$. 20 datasets were randomly generated for each setting and the mean and standard deviation are reported.}
    \label{tab:zero-sum-comparisons}
\end{table*}

Note that since the log-sum alternative $\min_{\bx} \{\Gamma(\bx)\}$ is a convex problem, both the Baseline and using gradient descent on top of our proposed approach to handle the exponential number of paths in Section \ref{sec:handle-exponential} are able to solve it to optimality. However, we see in Table \ref{tab:zero-sum-comparisons} the performance of the Baseline slightly tips off as the graph size $\lvert \cS \rvert$ increases, due to the fact that the objective in the baseline is estimated by sampling and the number of paths blows up exponentially with $\lvert \cS \rvert$.  

\begin{figure}[htp]
    \centering
    \includegraphics[width=0.7\textwidth]{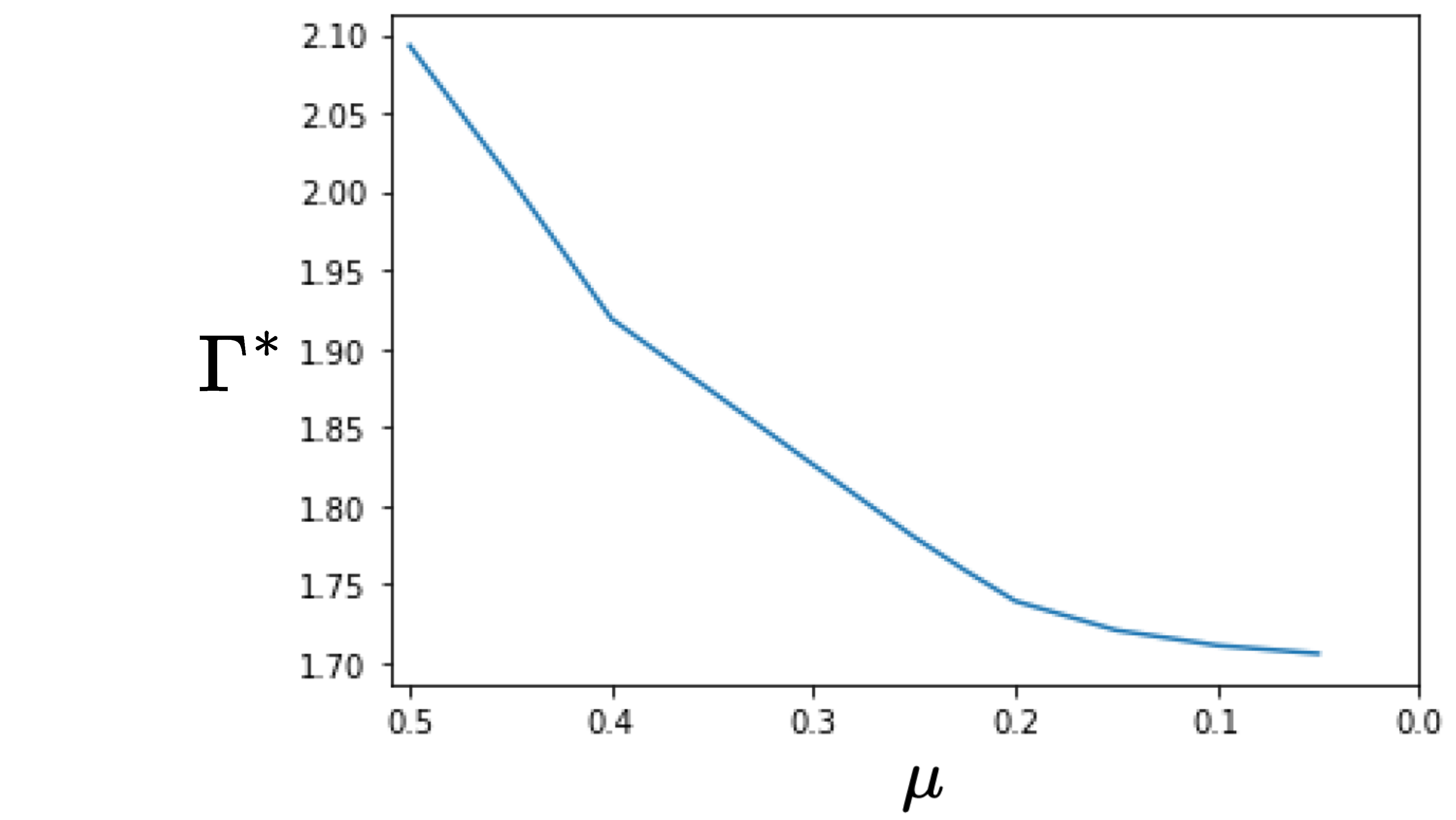}
    \caption{Optimal value  $\Gamma^*$  as a function of $\mu$ for a fixed synthetic dataset ($\lvert \cS \rvert$ = 50, p=0.8).}
    \label{fig:convergence}
\end{figure}

As the rationality of the adversary increases (associated with the decrease in $\mu$), we expect that the optimal defender reward will decrease as the adversary is able to take the best paths with a larger probability. Moreover, for a zero-sum game, by the guarantees in Proposition~\ref{prop:zero-sum} we can claim that the optimal solution would converge to the solution of \eqref{prob:zero-sum-shortest-path}. We note both the reward decrease and convergence in Figure~\ref{fig:convergence}.

\end{document}